\documentclass{amsart}

\usepackage{amsrefs}
\usepackage[all]{xy}
\usepackage{syntonly}
\usepackage{enumerate}
\usepackage{hyperref}
\usepackage{amsfonts}
\usepackage{amssymb}
\usepackage{amsmath}
\usepackage{amsthm}
\usepackage[inline]{enumitem}
\usepackage{tikz}
\usepackage{graphics}
\usepackage{graphicx}
\usepackage{color}
\usepackage{comment}

\usepackage{lineno}

\usepackage{selectp}

\newtheorem{theorem}{Theorem}
\newtheorem{corollary}[theorem]{Corollary}

\newtheorem{definition}[theorem]{Definition}
\newtheorem{lemma}[theorem]{Lemma}
\newtheorem{globalClaim}[theorem]{Claim}
\newtheorem{nonGlobalClaim}{Claim}[theorem]  

\newtheorem{observation}[theorem]{Observation}
\newtheorem{question}[theorem]{Question}

\newtheorem{fact}[theorem]{Fact}
\newtheorem{convention}[theorem]{Convention}
\newtheorem{remark}[theorem]{Remark}

\specialcomment{com}{ \color{red} }{\color{black}} 
\specialcomment{oldcom}{ \color{brown} }{\color{black}} 

\excludecomment{oldcom} 

\begin{document}
\title[Layered posets and Kunen's universal collapse]{Layered posets and Kunen's universal collapse}

\author{Sean D. Cox}
\email{scox9@vcu.edu}
\address{
Department of Mathematics and Applied Mathematics \\
Virginia Commonwealth University \\
1015 Floyd Avenue \\
Richmond, Virginia 23284, USA 
}

\thanks{The author thanks the anonymous referee for a careful review, and for the suggestion of adding some applications of the main iteration theorem.  The author also thanks Brent Cody and Monroe Eskew for helpful input about Section \ref{sec_Apps}.}


\begin{abstract}
We develop the theory of \emph{layered posets}, and  use the notion of layering to prove a new iteration theorem (Theorem \ref{thm_WeakCompact_Kunen_MainThm_Imprecise}): if $\kappa$ is weakly compact then \emph{\textbf{any}} universal Kunen iteration of $\kappa$-cc posets (each possibly of size $\kappa$) is $\kappa$-cc, as long as direct limits are used sufficiently often.  This iteration theorem simplifies and generalizes the various chain condition arguments for universal Kunen iterations in the literature on saturated ideals, especially in situations where finite support iterations are not possible.  We also  provide two applications: 
\begin{enumerate*}
 \item for any $n \ge 1$, a wide variety of $<\omega_{n-1}$-closed, $\omega_{n+1}$-cc posets of size $\omega_{n+1}$ can consistently be absorbed (as regular suborders) by quotients of saturated ideals on $\omega_n$ (see Theorem \ref{thm_IdealAbsorb} and Corollary \ref{cor_Examples}); and 
 \item for any $n \in \omega$, the Tree Property at $\omega_{n+3}$ is consistent with the Chang's Conjecture $(\omega_{n+3}, \omega_{n+1}) \twoheadrightarrow (\omega_{n+1}, \omega_n)$ (Theorem \ref{thm_ChangTree}).  
\end{enumerate*}
\end{abstract}

\maketitle

\tableofcontents

\section{Introduction}

A classic theorem of Solovay and Tennenbaum states:
\begin{theorem}[\cite{MR0294139}]\label{thm_FiniteSupport}
If $\kappa$ is regular uncountable, then any finite support iteration of $\kappa$-cc posets is $\kappa$-cc.
\end{theorem}

For iterations which are not finite support, the situation is much trickier.  A commonly used theorem in these more general situations is:
\begin{theorem}\label{thm_Jech1630}[see Theorem 16.30 of Jech~\cite{MR1940513}]
Let $\kappa$ be a regular uncountable cardinal and let $\alpha$ be a limit ordinal.  Let $\mathbb{P}_\alpha$ be an iteration such that for each $\beta < \alpha$, $\mathbb{P}_\beta = \mathbb{P}_\alpha \restriction \beta$ satisfies the $\kappa$-chain condition.  If $\mathbb{P}_\alpha$ is a direct limit and either
\begin{itemize}
 \item $\text{cf}(\alpha) \ne \kappa$ or
 \item $\text{cf}(\alpha) = \kappa$ and there are stationarily many $\beta < \alpha$ such that $\mathbb{P}_\beta$ is a direct limit
\end{itemize}
then $\mathbb{P}_\alpha$ satisfies the $\kappa$-chain condition.
\end{theorem}

The difficult part of applying Theorem \ref{thm_Jech1630} is typically the verification that $\mathbb{P}_\beta$ is $\kappa$-cc when $\mathbb{P}_\beta$ is an \emph{inverse} limit.  Usually this is taken care of simply by iterating \emph{small} posets, as in the following corollary, which is heavily used throughout the set theory literature:
\begin{corollary}\label{cor_IterateSmallPosets}[Baumgartner~\cite{MR823775}; see also Proposition 7.13 of Cummings~\cite{MR2768691}]
If $\kappa$ is inaccessible and $\langle \mathbb{P}_\alpha, \dot{Q}_\beta \ | \ \alpha \le \kappa, \beta < \kappa \rangle$ is an iteration such that:
\begin{itemize}
 \item For all $\alpha < \kappa$:  $\Vdash_{\mathbb{P}_\alpha} |\dot{Q}_\alpha| < \kappa$
 \item A direct limit is taken at $\kappa$ and on a stationary set of limit stages
below $\kappa$
\end{itemize}
Then $\mathbb{P}_\kappa$ has the $\kappa$-cc (in fact is $\kappa$-Knaster).
\end{corollary}

What if, in the statement of Corollary \ref{cor_IterateSmallPosets}, we don't assume that each $\dot{Q}_\alpha$ is forced to have size $<\kappa$?  This is in general a difficult problem, even if the posets $\dot{Q}_\beta$ are assumed to be highly closed; in fact there is a countable support iteration of length $\omega$ of countably closed, $\omega_2$-cc posets which fails to be $\omega_2$-cc (Exercise V.5.23 of Kunen~\cite{MR2905394}).  A prominent family of iterations which do \textbf{not} use small posets are the so-called \emph{universal Kunen iterations}.  Kunen introduced the first version in \cite{MR495118} to prove that saturated ideals on $\omega_1$ are consistent relative to huge cardinals.  There he assumed $\kappa$ was a huge cardinal and defined a $\kappa$-length finite support iteration $\langle \mathbb{P}_\alpha, \dot{Q}_\beta \ | \ \alpha \le \kappa, \beta < \kappa \rangle$ where, importantly, each $\dot{Q}_\alpha$ was forced to have size $\kappa$ and was in fact from an inner model of $V^{\mathbb{P}_\alpha}$.\footnote{Namely, $\dot{Q}_\alpha$ is the Silver collapse to turn $\kappa$ into $\alpha^+$, as defined in the model $V^{V_\alpha \cap \mathbb{P}_\alpha}$ in the case that $V_\alpha \cap \mathbb{P}_\alpha$ is a regular suborder of $\mathbb{P}_\alpha$, and trivial otherwise.  See Section \ref{sec_KunenUniversal} for a discussion of such iterations.}  He used the Knaster property of the Silver collapse, together with Theorem \ref{thm_FiniteSupport} about finite support iterations, to prove that $\mathbb{P}_\kappa$ had the $\kappa$-cc.  The poset $\mathbb{P}_\kappa$ is highly universal, even for many posets of size $\kappa$; contrast this with the Levy Collapse $\text{Col}(\omega, < \kappa)$ which is only universal for posets of size $<\kappa$.  The strong universality property of $\mathbb{P}_\kappa$ is the key to obtaining master conditions which enable one to construct a saturated ideal in the final model.

Laver generalized Kunen's universal iteration to contexts where larger supports were used, in order to obtain saturated ideals on $\omega_2$ and beyond (see \cite{MattHandbook} for a discussion of the history).  However in those settings the $\kappa$-cc preservation of the resulting universal iteration had to be checked on a case-by-case basis,\footnote{e.g.\ Lemma 4.8 of \cite{MR2151585}; \cite{MR662045} } and seemingly was believed to be a somewhat delicate issue (see Remark 20.3 of Cummings~\cite{MR2768691}).   We prove that, at least if $\kappa$ is weakly compact, then the situation isn't really delicate at all; \textbf{any} universal Kunen iteration of $\kappa$-cc posets---that is, where $\dot{Q}_\beta$ is forced by $V_\beta \cap \mathbb{P}_\beta$ to be $\kappa$-cc\footnote{Note this is an apparently weaker assumption than requiring that $\dot{Q}_\beta$ is forced by $\mathbb{P}_\beta$ to be $\kappa$-cc, though in the end they are equivalent if $\kappa$ is weakly compact. }---will be $\kappa$-cc, provided that direct limits are used often.  

The proof heavily uses the notion of a \emph{layered poset}.  Layering has previously appeared in the literature in the form of so-called \emph{layered ideals} (e.g.\ \cite{MR942519}), but it is a general property of posets which is quite useful and interesting.  A poset $\mathbb{P}$ is $\kappa$ stationarily layered iff there are stationarily many $\mathbb{Q} \in P_\kappa(\mathbb{P})$ which are regular suborders of $\mathbb{P}$; $\kappa$ club layering is defined similarly.

Among the basic facts we prove about layering is the following lemma, proved in Section \ref{sec_RelationToKnaster}:
\begin{lemma}\label{lem_LayeredImpliesKnaster}
If $\mathbb{P}$ is $\kappa$ stationarily layered then $\mathbb{P}$ is $\kappa$-Knaster.
\end{lemma}
The converse is, in general, false; see Cox-L\"ucke~\cite{Cox_Luecke}.  In \cite{Cox_Luecke} it is also proved that weak compactness of $\kappa$ is equivalent to ``every $\kappa$-cc poset is $\kappa$ stationarily layered".  In short:
\begin{equation*}
\kappa \text{-stationarily layered } \implies \kappa\text{-Knaster } \implies \kappa \text{-cc}
\end{equation*}
and if $\kappa$ is weakly compact then the 3 notions are equivalent (and highly productive).

In this paper, the notion of layering is used to prove the following general theorems about universal Kunen iterations.  See Section \ref{sec_KunenUniversal} for a general discussion of universal Kunen iterations, and precise statements of the theorems.

\begin{theorem}\label{thm_MainTheorem_Mahlo_Imprecise}
(See Theorem \ref{thm_MainTheorem_Mahlo_Precise} for the precise statement)  If $\kappa$ is Mahlo then any universal Kunen iteration of sufficiently $\kappa$-layered posets is $\kappa$-stationarily layered, provided direct limits were taken at all inaccessible $\gamma \le \kappa$.  
\end{theorem}

\begin{theorem}\label{thm_WeakCompact_Kunen_MainThm_Imprecise}
(See Corollary \ref{thm_WeakCompact_Kunen_MainThm_Precise} for the precise statement)  If $\kappa$ is weakly compact then any universal Kunen iteration of $\kappa$-cc posets is $\kappa$-stationarily layered, provided direct limits were taken at all inaccessible $\gamma \le \kappa$.  
\end{theorem}

Note that, in particular, by Lemma \ref{lem_LayeredImpliesKnaster} the universal Kunen iteration is $\kappa$-Knaster in the conclusion of both theorems.

We give several applications of Theorem \ref{thm_WeakCompact_Kunen_MainThm_Imprecise}.  The statement of the meta-mathematical Theorem \ref{thm_IdealAbsorb} below is necessarily somewhat technical; the role of the parameter $r$ in \eqref{eq_DefinablePlus2} is to allow, for example, that the forcing is defined inside $H_{\kappa^{++}}$ using some bookkeeping device; e.g. a parameter which is a $\kappa^+$-to-one surjection from $\kappa^+ \to H_{\kappa^+}$ used to define some iteration of length $\kappa^+$.  \textbf{See Corollary \ref{cor_Examples} for specific examples of posets which satisfy the requirements of Theorem \ref{thm_IdealAbsorb}.}
\begin{theorem}\label{thm_IdealAbsorb}
Assume $\phi(-,-)$ is a set theoretic formula with two free variables, and ZFC proves the following:  whenever $\kappa$ is a successor of a regular cardinal then there exists an $r \in H_{\kappa^{++}}$ such that, letting $\mu$ be the cardinal predecessor of $\kappa$:
\begin{align}\label{eq_DefinablePlus2}
\begin{split}
 \mathbb{S}^{H_{\kappa^{++}}}_{r,\phi}:= \big\{ z \in H_{\kappa^{++}} \ |   (H_{\kappa^{++}}, \in) \models \phi(z,r) \big\}  
\text{ is a } <\mu \text{-closed, } \\ 
 \kappa^+ \text{-cc poset} \text{ of size at most } \kappa^+
\end{split}
\end{align}

\noindent \textbf{Then:}  whenever $V$ is a model of ZFC and $V \models $ ``$\mu < \kappa$ are regular cardinals and $\kappa$ is a huge cardinal", there is a $<\mu$-closed poset in $V$ which forces the following statements:
\begin{itemize}
 \item $\kappa = \mu^+ = 2^\mu$ and $2^\kappa = \kappa^+$;
 \item There is a normal, saturated ideal $\mathcal{J}$ on $\kappa$;
 \item There exists some $r \in H_{\kappa^{++}}$ such that \eqref{eq_DefinablePlus2} holds, and there exists some regular embedding
\[
e: \mathbb{S}^{H_{\kappa^{++}}}_{r,\phi} \to \wp(\kappa)/\mathcal{J}
\]
\end{itemize}
\end{theorem}

The formula $\phi(-,-)$ in the assumptions of Theorem \ref{thm_IdealAbsorb} can be designed in ways to specify a particular $\omega_n$, or to specify that certain other requirements are met, if desired.\footnote{Note that $\kappa$ is definable in $(H_{\kappa^{++}}, \in)$, as ``the 2nd largest cardinal".  So, for example, if one has in mind a definable $\kappa^+$-cc forcing which only has the desired properties when, say, $\kappa = \omega_2$, then one can simply insert ``and the 2nd largest cardinal is $\omega_2$" into the formula $\phi(-,-)$.  This will have the effect of making $\mathbb{S}^{H_{\kappa^{++}}}_{r,\phi}$ into the trivial poset when $\kappa \ne \omega_2$.  Similarly, if one has in mind a definable $\kappa^+$-cc forcing which only makes sense when $2^\kappa= \kappa^+$, then simply insert ``and the 2nd largest cardinal satisfies GCH" into the formula $\phi$; this would have the effect of making $\mathbb{S}^{H_{\kappa^{++}}}_{r,\phi}$ into the trivial poset when $2^\kappa > \kappa^+$. }  The following corollary lists some specific examples which satisfy the requirements of Theorem \ref{thm_IdealAbsorb}.  We say that $\mathbb{R}$ absorbs $\mathbb{Q}$ if there exists a regular embedding from $\mathbb{Q} \to \mathbb{R}$.  

\begin{corollary}\label{cor_Examples}
Fix $n \in \omega$, $n \ge 1$.  Each of the following posets can consistently be absorbed by the quotient of a saturated ideal on $\omega_{n}$:
\begin{enumerate}
 \item Jensen's $<\omega_n$-closed forcing of size $\omega_{n+1}$ (under GCH) to add an $\omega_{n}$-Kurepa tree (see \cite{MR2768691}).
 \item The forcing to add a single Hechler subset of $\omega_{n}$ (\cite{MR1355135}).
 \item A $<\omega_{n-1}$-closed, $\omega_{n+1}$-Suslin tree.
 \item The $\omega_{n+1}$-length, $<\omega_{n-1}$-support iteration of adding Hechler subsets of $\omega_{n-1}$ (\cite{MR1355135}).  More generally, any sufficiently definable\footnote{See Theorem \ref{thm_IdealAbsorb} for precisely what is meant here.} $<\omega_{n-1}$-support iteration of length $\omega_{n+1}$ which uses $<\omega_{n-1}$-closed posets of size $\le \omega_{n}$ at each step. 
 \item $\le \omega_{n-1}$-support, $\omega_{n+1}$-length iteration of $\text{Sacks}(\omega_{n-1})$ (as in Kanamori~\cite{MR593029})  
 \item (the following poset can consistently be absorbed into the quotient of a saturated ideal on $\omega_2$):  The $\sigma$-closed, $\omega_2$-cc poset of size $\omega_3$ from Shelah~\cite{Shelah_WeakGenMA}, which he used to force $2^{\omega_1} = \omega_3$ together with a generalized version of Martin's Axiom for a certain subclass of the $\sigma$-closed, $\omega_2$-cc posets. 
\end{enumerate}
\end{corollary}

Finally we prove that a certain instance of Chang's Conjecture is consistent the Tree Property:
\begin{theorem}\label{thm_ChangTree}
If $\mu < \kappa$ are both regular and $\kappa$ is a huge cardinal, then there is a forcing extension which satisfies the following:
\begin{itemize}
 \item $\kappa = \mu^+$;
 \item the Chang's Conjecture $(\mu^{+3}, \mu^+) \twoheadrightarrow (\mu^+, \mu)$ holds; and
 \item The Tree Property holds at $\mu^{+3}$.
\end{itemize}
In particular, for any $n \in \omega$, the tree property at $\omega_{n+3}$ is consistent with the Chang's Conjecture
\[
(\omega_{n+3}, \omega_{n+1}) \twoheadrightarrow (\omega_{n+1}, \omega_n)
\] 
\end{theorem}

The paper is organized as follows:  
\begin{itemize}
 \item Section \ref{sec_Prelims} includes some preliminaries;
 \item Section \ref{sec_LayeredPosets} introduces layered posets; proves that layering implies the Knaster condition; provides some examples; and proves a useful lemma about layering and elementary embeddings. 
 \item Section \ref{sec_InductReductions} introduces the important notion of a \emph{coherent, conservative} system of reduction operations and proves Lemma \ref{lem_WhatToDoAtLimits}, which is crucial to passing limit stages in the proofs of Section \ref{sec_KunenUniversal};
 \item Section \ref{sec_KunenUniversal} provides the relevant background for universal Kunen-type iterations, and proves (the precise versions of) Theorems \ref{thm_MainTheorem_Mahlo_Imprecise} and \ref{thm_WeakCompact_Kunen_MainThm_Imprecise}.
 \item Section \ref{sec_Apps} proves Theorems \ref{thm_IdealAbsorb} and \ref{thm_ChangTree}.
\end{itemize}

\section{Preliminaries}\label{sec_Prelims}

\subsection{Regular embeddings and quotient forcing}

If $p,q$ are conditions in a poset $\mathbb{P}$, we write $p \le q$ to mean that $p$ is stronger than $q$, $p \parallel q$ to mean that $p$ and $q$ are compatible, and $p \perp q$ to mean that $p$ and $q$ are incompatible.  Given $\mathbb{P}$-names $\tau$ and $\tau'$, we will say that $\tau$ is equivalent to $\tau'$ (modulo $\mathbb{P}$) iff $\Vdash_{\mathbb{P}} \tau = \tau'$.  The \emph{Maximality Principle} is the fact that if $p \Vdash_{\mathbb{P}} \exists x \ \phi(x)$ where $\phi$ is a statement in the forcing language, then there is some $\mathbb{P}$-name $\tau$ such that $p \Vdash \phi(\tau)$.

A function $e: \mathbb{P} \to \mathbb{Q}$ is a \emph{regular embedding} iff it is order and incompatibility preserving, and whenever $A$ is a maximal antichain in $\mathbb{P}$ then $e[A]$ is a maximal antichain in $\mathbb{Q}$.  If $\text{id}: \mathbb{P} \to \mathbb{Q}$ is a regular embedding then $\mathbb{P}$ is called a \emph{regular suborder} of $\mathbb{Q}$.  If $\text{id}:\mathbb{P} \to \mathbb{Q}$ is order and incompatibility preserving---but not necessarily regular---then we say $\mathbb{P}$ is a suborder of $\mathbb{Q}$.

\begin{definition}
Let $\mathbb{P}$ be a suborder of $\mathbb{Q}$ and $q \in \mathbb{Q}$.  We say that $p \in \mathbb{P}$ is a \emph{reduct of $q$ into $\mathbb{P}$} iff   
\begin{equation*}
  \forall p' \le_{\mathbb{P}} p \ p' \parallel_{\mathbb{Q}} q.
\end{equation*}
If $p$ is also $\ge_{\mathbb{Q}} q$ then we say $p$ is a \emph{nice} reduct of $q$ into $\mathbb{P}$.
\end{definition}

The following lemma is standard, but we provide the short proof for convenience.  Note that the reduction characterization yields a $\Sigma_0$ characterization of the regular suborder relation.
\begin{lemma}\label{lem_ReductCharactOfReg}
Let $\mathbb{P}$ be a suborder of $\mathbb{Q}$.  The following are equivalent:
\begin{enumerate}
  \item $\mathbb{P}$ is a regular suborder of $\mathbb{Q}$;
  \item For every $q \in \mathbb{Q}$ there is a reduct of $q$ into $\mathbb{P}$
\end{enumerate}
\end{lemma}
\begin{proof}
Assume $\mathbb{P}$ is a regular suborder of $\mathbb{Q}$.  So every maximal antichain in $\mathbb{P}$ is a maximal antichain in $\mathbb{Q}$.  Assume for a contradiction there is some $q \in \mathbb{Q}$ which has no reduct; this implies that
\begin{equation*}
D_{\perp q}:= \{ p \in \mathbb{P} \ | \  p \perp_{\mathbb{Q}} q \} \text{ is dense in } \mathbb{P}.
\end{equation*}
Let $G$ be $\mathbb{Q}$-generic with $q \in G$.  Then $G \cap \mathbb{P}$ is $\mathbb{P}$-generic and so $G \cap D_{\perp q} \ne \emptyset$, contradicting that $G$ is a filter on $\mathbb{Q}$ and that $q \in G$.

For the other direction, let $A$ be a maximal antichain in $\mathbb{P}$.  Let $q$ be any condition in $\mathbb{Q}$, and let $p_q$ be any reduct of $q$ into $\mathbb{P}$.  Let $a \in A$ be compatible in $\mathbb{P}$ with $p_q$, and let $p' \in \mathbb{P}$ witness this compatibility.  Because $p' \le p_q$, it follows that $p' \parallel_{\mathbb{Q}} q$.  Then $q$ is compatible with $a \in A$.  Because $q$ was arbitrary then $A$ is maximal in $\mathbb{Q}$.
\end{proof}

\begin{definition}\label{def_GeneralizedOrder}
Suppose $e: \mathbb{P} \to \mathbb{Q}$ is a regular embedding, $p \in \mathbb{P}$, and $q \in \mathbb{Q}$.  We write $q \le_e p$ to mean:
\begin{equation*}
q \Vdash_{\mathbb{Q}} \check{p} \in \Big(\check{e}^{-1}[\dot{G}_{\mathbb{Q}}]\Big) \uparrow
\end{equation*}
where $\Big( \check{e}^{-1}[\dot{G}_{\mathbb{Q}}]\Big) \uparrow$ denotes the set of all conditions in $\mathbb{P}$ which are weaker than some condition in $\check{e}^{-1}[\dot{G}_{\mathbb{Q}}]$.

If $q \in \mathbb{Q}$ and $p \in \mathbb{P}$, we say that $p$ is an \emph{$e$-reduct of $q$} iff
\[
\forall p' \le_{\mathbb{P}} p \ \ e(p') \parallel_{\mathbb{Q}} q  .
\]
\end{definition}
Note that if $e = \text{id}$ and $\mathbb{Q}$ is separative, then $q \le_e p$ is equivalent to $q \le_{\mathbb{Q}} p$.

\begin{lemma}\label{lem_GeneralizedNiceReduct}
Suppose $\mathbb{B}$ and $\mathbb{C}$ are complete boolean algebras, and $e: \mathbb{B} \to \mathbb{C}$ is a regular embedding.  Let $c \in \mathbb{C}$ and set
\begin{equation*}
b^e_c:= \text{sup}_{\mathbb{B}} \{ b \in \mathbb{B} \ | \ b \text{ is an } e \text{-reduct of } c   \} .
\end{equation*}
Then
\begin{enumerate}
 \item $b^e_c$ is an $e$-reduct of $c$ and moreover is the largest one;
 \item $b^e_c \ge_e c$ (as in Definition \ref{def_GeneralizedOrder})
\end{enumerate}  
\end{lemma}
\begin{proof}
To see that $b^e_c$ is an $e$-reduct of $c$, let $b' \le b^e_c$ with $b' \ne 0$.  By the definition of $b^e_c$ there is some $b$ which is an $e$-reduct of $c$ such that $b' \wedge b \ne 0$.  Then $b' \wedge b \le b$ and because $b$ is an $e$-reduct of $c$ we have $e(b' \wedge b) \parallel c$.  Because $e$ is order preserving and $b' \wedge b \le b'$ we have $e(b') \parallel c$, which completes the proof that $b^e_c$ is an $e$-reduct of $c$.  Clearly it is the largest such, by its definition as the sup of all $e$-reducts of $c$.

To see that $b^e_c \le_e c$: suppose toward a contradiction that this were false; then there is some $c' \le c$ such that
\begin{equation*}
c' \Vdash_{\mathbb{C}} b^e_c \notin e^{-1}[\dot{G}_{\mathbb{C}}] .
\end{equation*}
It follows easily that
\begin{equation}\label{eq_cprimePerp}
c' \perp e(b^e_c).
\end{equation}
Let $b'$ be any $e$-reduct of $c'$.  Then $b' \vee b^e_c > b^e_c$; otherwise $b' \le b^e_c$ and because $c' \parallel e(b')$ this would contradict \eqref{eq_cprimePerp}.  Because $b' \vee b^e_c > b^e_c$ and $b^e_c$ is defined as the supremum of all $e$-reducts of $c$, then $b' \vee b^e_c$ is \textbf{not} a reduct of $c$.  We will prove that $b' \vee b^e_c$ is in fact an $e$-reduct of $c$ which will yield a contradiction and complete the proof.  So let $r \le b' \vee b^e_c$, $r \ne 0$.  Then $r \parallel b'$ or $r \parallel b^e_c$.  
\begin{itemize}
 \item CASE 1:  $r \parallel b'$.  Because $b'$ is an $e$-reduct of $c'$ then $e(r) \parallel c'$, and because $c' \le c$ this implies $e(r) \parallel c$.
 \item CASE 2:  $r \parallel b^e_c$.  Then because  $b^e_c$ is an $e$-reduct of $c$ (by the first part of the proof), we have $e(r) \parallel c$.
\end{itemize}
\end{proof}

\begin{corollary}\label{cor_b_c_is_nicereduct}
Suppose $\mathbb{B}$ is a complete regular subalgebra of $\mathbb{C}$, and $c \in \mathbb{C}$.  Let 
\begin{equation*}
b_c:= \text{sup}_{\mathbb{B}} \{  b \in \mathbb{B} \ | \ b \text{ is a reduct of } c \}.
\end{equation*}
Then $b_c$ is a nice reduct of $c$.
\end{corollary}
\begin{proof}
This follows immediately from Lemma \ref{lem_GeneralizedNiceReduct} by taking the $e$ from that lemma to be the identity map from $\mathbb{B} \to \mathbb{C}$. 
\end{proof}

The following fact is well-known:
\begin{lemma}\label{lem_LiftRegular}
If $e: \mathbb{P} \to \mathbb{Q}$ is a regular embedding of separative posets, then $e$ lifts to a regular embedding of their Boolean completions.
\end{lemma}

If the map $e$ from Lemma \ref{lem_LiftRegular} is the identity---i.e.\ if $\mathbb{P}$ is a regular suborder of $\mathbb{Q}$---then we can essentially take the lifting to the boolean completions to be the identity map, which simplifies the notation a bit.  In other words, if $\mathbb{P}$ is a regular suborder of $\mathbb{Q}$ then we can extend both to view $\text{ro}(\mathbb{P})$ as a regular suborder of $\text{ro}(\mathbb{Q})$.  I'm grateful to Joel David Hamkins for pointing this out to me on MathOverflow:
\begin{lemma}[\cite{206871}]\label{lem_HamkinsLemma}
Suppose $\mathbb{P}$ is a regular suborder of the separative poset $\mathbb{Q}$.  In $\text{ro}(\mathbb{Q})$ define
\begin{equation*}
\mathbb{B}(\mathbb{P}):= \big\{ \text{sup}_{\text{ro}(\mathbb{Q})}(X) \ | \ X \subseteq \mathbb{P}  \big\}.
\end{equation*}
Then $\mathbb{B}(\mathbb{P})$ is a complete regular subalgebra of $\text{ro}(\mathbb{Q})$, and $\mathbb{P}$ is dense in $\mathbb{B}(\mathbb{P})$ (so $\mathbb{B}(\mathbb{P}) \simeq \text{ro}(\mathbb{P})$). 
\end{lemma}

\begin{definition}\label{def_WeakestReduct}
Let $\mathbb{P}$ be a regular suborder of the separative poset $\mathbb{Q}$, and let $q \in \mathbb{Q}$.  Set 
 \begin{equation*}
 \text{WeakestReduct}_{\text{ro}(\mathbb{P})}(q) = \text{sup} \{ b \in \text{ro}(\mathbb{P}) \ | \ b \text{ is a reduct of } q   \}
 \end{equation*}
where we are viewing $\text{ro}(\mathbb{P})$ as a regular complete subalgebra of $\text{ro}(\mathbb{Q})$ as in Lemma \ref{lem_HamkinsLemma}.
\end{definition}

We will often use the following convention:
\begin{convention}\label{conv_ViewNames}
Assume $\mathbb{P}$ is a regular suborder of $\mathbb{Q}$, $\sigma$ is a statement in the $\mathbb{P}$ forcing language, and $q \in \mathbb{Q}$.  We will write $q \Vdash_{\mathbb{Q}} \sigma$ to mean the following: whenever $G$ is $\mathbb{Q}$-generic with $q \in G$, then $V[G \cap \mathbb{P}] \models \sigma_{G \cap \mathbb{P}}$. 
\end{convention}

\begin{lemma}\label{lem_IfWeakestReductForcesThenOrigDoes}
Let $\mathbb{P}$ be a regular suborder of $\mathbb{Q}$, $q \in \mathbb{Q}$, and $\sigma$ a statement in the forcing language of $\mathbb{P}$.  Then the following are equivalent:
\begin{enumerate}
 \item\label{item_q_forces_sigma} $q \Vdash_{\mathbb{Q}} \sigma$ (as in Convention \ref{conv_ViewNames});
 \item\label{item_WeakestReduct_forces_sigma} $\text{WeakestReduct}_{\text{ro}(\mathbb{P})}(q) \Vdash_{\text{ro}(\mathbb{P})} \sigma$
\end{enumerate}
\end{lemma}
\begin{proof}
By Corollary \ref{cor_b_c_is_nicereduct}, $\text{WeakestReduct}_{\text{ro}(\mathbb{P})}(q)$ is a nice reduct of $q$ into $\text{ro}(\mathbb{P})$; i.e.\ it is a reduct and moreover $\text{WeakestReduct}_{\text{ro}(\mathbb{P})}(q) \ge q$.  

To see (\ref{item_q_forces_sigma}) implies (\ref{item_WeakestReduct_forces_sigma}):  if not then there is some $p' \le \text{WeakestReduct}_{\text{ro}(\mathbb{P})}(q)$ in $\text{ro}(\mathbb{P})$ such that $p'$ forces the negation of $\sigma$.  But because $\text{WeakestReduct}_{\text{ro}(\mathbb{P})}(q)$ is a reduct then $p'$ is compatible in $\text{ro}(\mathbb{Q})$ with $q$, yielding a contradiction to assumption (\ref{item_q_forces_sigma}).

To see that (\ref{item_WeakestReduct_forces_sigma}) implies (\ref{item_q_forces_sigma}):  this is simply because $\text{WeakestReduct}_{\text{ro}(\mathbb{P})}(q) \ge q$; any $\text{ro}(\mathbb{Q})$-generic which includes $q$ must also include $\text{WeakestReduct}_{\text{ro}(\mathbb{P})}(q)$, so assumption (\ref{item_WeakestReduct_forces_sigma}) ensures that $\sigma$ holds in any such generic extension.
\end{proof}

If $\mathbb{P}$ is a regular suborder of $\mathbb{Q}$ and $G_{\mathbb{P}}$ is generic for $\mathbb{P}$, then in $V[G_{\mathbb{P}}]$ the quotient $\mathbb{Q}/G_{\mathbb{P}}$ is defined as the set of all $q \in \mathbb{Q}$ which are compatible (in $\mathbb{Q}$) with every member of $G_{\mathbb{P}}$; the ordering of $\mathbb{Q}/G_{\mathbb{P}}$ is the order inherited from $\mathbb{Q}$.  The following fact is standard:
\begin{fact}\label{fact_OneStepTwoSteps}
If $\mathbb{P}$ is a regular suborder of $\mathbb{Q}$, then 
\[
\mathbb{Q} \text{ is forcing equivalent to } \mathbb{P}  \ * \ \check{\mathbb{Q}}/\dot{G}_{\mathbb{P}}.
\]
\end{fact}

We will also use the following fact in the proof of Theorem  \ref{thm_IdealAbsorb}:
\begin{fact}\label{fact_CannotAddAntichains}
If $\mathbb{P}$ is a regular suborder of $\mathbb{Q}$ and $\mathbb{Q}$ has the $\kappa$-cc, then 
\[
\Vdash_{\mathbb{P}} \ \ \check{\mathbb{Q}}/\dot{G}_{\mathbb{P}} \text{ has the } \kappa \text{-cc}.
\]
\end{fact}
\begin{proof}
Suppose for a contradiction that there is some $p \in \mathbb{P}$ forcing that there is a $\kappa$-sized antichain in the quotient; by the Maximaility Principle there is a $\mathbb{P}$-name $\langle \dot{q}_\xi \ | \ \xi < \kappa \rangle$ which is forced by $p$ to be an antichain in the quotient.  But then $\{ (p,\dot{q}_\xi) \ | \ \xi < \kappa \}$ is easily seen to be an antichain in $\mathbb{P} * \check{\mathbb{Q}}/\dot{G}_{\mathbb{P}}$.  By Fact \ref{fact_OneStepTwoSteps} this implies $\mathbb{Q}$ fails to have the $\kappa$-cc, a contradiction.
\end{proof}

\subsection{Generalized Sacks forcing}

The proof of Theorem \ref{thm_ChangTree} will make use of the following generalized Sacks forcing.  If $T$ is a subtree of ${}^{<\kappa} 2$, $\alpha < \kappa$, and $s \in {}^\alpha 2$, we say that \emph{$s$ splits in $T$} iff both $s^\frown 0$ and $s^\frown 1$ are in $T$.  
\begin{definition}[Kanamori~\cite{MR593029}, Definition 1.1]\label{def_Sacks}
Let $\kappa$ be a regular uncountable cardinal.  $\text{Sacks}(\kappa)$ is the poset of all subtrees $T$ of ${}^{<\kappa} 2$ such that:
\begin{enumerate}
 \item\label{item_Closed} If $\alpha < \kappa$ is a limit ordinal, $s \in {}^\alpha 2$, and $s \restriction \beta \in T$ for every $\beta < \alpha$, then $s \in T$.
 \item\label{item_Splits} If $s \in T$ then there is a $t \in T$ such that $s \subset t$ and $t$ splits in $T$.
 \item\label{item_ClubSplitLevels} If $\alpha < \kappa$ is a limit ordinal, $s \in {}^\alpha 2$, and $s \restriction \beta$ splits in $T$ for cofinally many $\beta <\alpha$, then $s$ splits in $T$.
\end{enumerate}
$T$ is stronger than $S$ iff $T \subseteq S$.
\end{definition}

\begin{lemma}\label{lem_ClubSplittingBranch}
If $T \in \text{Sacks}(\kappa)$ and $f: \kappa \to 2$ is a cofinal branch through $T$ (i.e. $f \restriction \alpha \in T$ for all $\alpha < \kappa$), then 
\[
C_{f,T}:= \{ \alpha < \kappa \ | \  f \restriction \alpha \text{ splits in } T  \} \text{ is a closed, unbounded subset of } \kappa.
\]
\end{lemma}
\begin{proof}
By requirement \ref{item_ClubSplitLevels} of Definition \ref{def_Sacks} it suffices to show that $C_{f,T}$ is cofinal in $\kappa$.  Fix any $\alpha_0 < \kappa$.  Because $f \restriction \alpha_0 \in T$ then by requirement \ref{item_Splits} of Definition \ref{def_Sacks} there is some $g \supseteq f \restriction \alpha_0$ which splits in $T$; i.e.\ $g^\frown 0 \in T$ and $g^\frown 1 \in T$.  Let $\beta_{f,g}$ be the maximal ordinal such that $g \restriction \beta_{f,g} = f \restriction \beta_{f,g}$.  Notice that
\[
\beta_{f,g} \ge \alpha_0.
\]
We claim that $f \restriction \beta_{f,g}$ splits in $T$.  For conceptual clarity, consider two cases:  
\begin{enumerate*}
 \item If $\beta_{f,g} = \text{dom}(g)$ then $
 g^\frown 0 = (f \restriction \beta_{f,g})^\frown 0$ and $g^\frown 1 = (f \restriction \beta_{f,g})^\frown 1$.  Because both of these functions are in $T$ by the assumption that $g$ splits in $T$, this shows that $f \restriction \beta_{f,g}$ splits in $T$.
 \item If $\beta_{f,g} < \text{dom}(g)$, then $g \restriction \beta_{f,g} = f \restriction \beta_{f,g}$ and $g(\beta_{f,g}) \ne f(\beta_{f,g})$.  Because $g \restriction (\beta_{f,g} + 1)$ and $f \restriction (\beta_{f,g} + 1)$ are both in $T$ this again implies that $f \restriction \beta_{f,g}$ splits in $T$.
\end{enumerate*}

\end{proof}

Kanamori (\cite{MR593029}, Lemma 1.2) proves that $\text{Sacks}(\kappa)$ is $<\kappa$ closed, but we will actually need a bit more.  Recall that if $D$ is a subset of a poset $\mathbb{P}$, we say $D$ is \emph{directed} iff whenever $d_0,d_1 \in D$ then there is a $d \in D$ such that $d \le d_0$ and $d \le d_1$.  A poset $\mathbb{P}$ is $<\kappa$-directed closed iff whenever $D$ is a $<\kappa$-sized, directed subset of $\mathbb{P}$, then $D$ has a lower bound in $\mathbb{P}$.

\begin{lemma}\label{lem_SacksDirectedClosed}
For all regular uncountable $\kappa$, the poset $\text{Sacks}(\kappa)$ is $<\kappa$-directed closed.
\end{lemma}
\begin{proof}
Let $D$ be a $<\kappa$-sized, directed subset of $\text{Sacks}(\kappa)$.  We prove that $\bigcap D$ inherits all requirements of conditionhood from Definition \ref{def_Sacks}.  The only nontrivial part is verifying requirement \ref{item_Splits}; all other requirements of Definition \ref{def_Sacks} are easily inherited by arbitrary intersections of conditions.

So suppose $s \in \bigcap D$; we need to find some $t \supset s$ which splits in $\bigcap D$.  
\begin{nonGlobalClaim}
There is an $f: \kappa \to 2$ extending $s$ such that $f$ is a cofinal branch through every member of $D$.
\end{nonGlobalClaim}  
\begin{proof}
We recursively define $f \restriction \alpha$ for all $\alpha \in [\text{dom}(s), \kappa]$ and verify inductively that $f \restriction \alpha \in \bigcap D$ for all such $\alpha$.  First set $f \restriction \text{dom}(s):= s$.  If $\alpha$ is a limit ordinal then $f \restriction \alpha$ is defined as the union of $f \restriction \beta$ for all $\beta < \alpha$; the induction hypothesis and requirement \ref{item_Closed} of Definition \ref{def_Sacks} ensure that $f \restriction \alpha \in \bigcap D$.  

Now suppose $f \restriction \alpha$ is defined; we want to define $f \restriction (\alpha + 1)$.  We claim that at least one of $(f \restriction \alpha)^\frown 0$ or $(f \restriction \alpha)^\frown 1$ is an element of $\bigcap D$.  If not, then there are $T_0, T_1 \in D$ such that $(f \restriction \alpha)^\frown 0 \notin T_0$ and $(f \restriction \alpha)^\frown 1 \notin T_1$.  Because $D$ is directed, there is a $T \in D$ such that $T \subseteq T_0 \cap T_1$.  Now the induction hypothesis ensures that $f \restriction \alpha \in T$; and requirement \ref{item_Splits} of Definition \ref{def_Sacks} ensures that $(f \restriction \alpha)^\frown i \in T$ for at least one $i \in \{ 0,1 \}$.  But then $(f \restriction \alpha)^\frown i \in T_i$, a contradiction.
\end{proof}

Let $f: \kappa \to 2$ be as given by the claim, and for every $T \in D$ let $C_{f,T}$ be the club from Lemma \ref{lem_ClubSplittingBranch}.  Because $|D|<\kappa$, $C:= \bigcap_{T \in D} C_{f,T}$ is club in $\kappa$.  Pick some $\beta \in C \cap (\text{dom}(s),\kappa)$.  Then $f \restriction \beta$ splits in every member of $D$; i.e.
\[
(f \restriction \beta)^\frown 0 \text{ and } (f \restriction \beta)^\frown 1 \text{ are both in } \bigcap D.
\]
Because $f \restriction \beta \supset s$ this completes the proof.
\end{proof}

\subsection{Other miscellaneous background}

A regular uncountable cardinal $\kappa$ is:
\begin{itemize}
 \item \emph{weakly compact} iff $\kappa$ is inaccessible and whenever $M$ is a transitive $ZF^-$ model of size $\kappa$ and $M$ is closed under $<\kappa$ sequences, then there is a transitive set $N$ and an elementary embedding $j: M \to N$ with critical point $\kappa$.  There are many other well-known equivalent formulations.
 \item \emph{almost huge} iff there is an elementary embedding $j: V \to N$ with critical point $\kappa$ such that $N$ is closed under $<j(\kappa)$ sequences.
 \item \emph{huge} iff there is an elementary embedding $j: V \to N$ with critical point $\kappa$ such that $N$ is closed under $j(\kappa)$ sequences.
\end{itemize}

The \emph{Tree Property} holds at $\kappa$ iff whenever $T$ is a tree of height $\kappa$ and each level of $T$ is of size $<\kappa$, then $T$ has a cofinal branch.  If 
\[
\mu_3 > \mu_2 > \mu_1 > \mu_0
\]
are cardinals then the ``Chang's Conjecture" $(\mu_3, \mu_2) \twoheadrightarrow (\mu_1, \mu_0)$ means:  whenever $\mathfrak{A}$ is a first order structure on $\mu_3$ in a countable language, then there is an $X \prec \mathfrak{A}$ such that $|X| = \mu_1$ and $|X \cap \mu_2| = \mu_0$.  Various versions of this are discussed at length in \cite{MattHandbook}.

The proof of Theorem \ref{thm_IdealAbsorb} makes use of Foreman's Duality Theorem, specifically regarding the ``Magidor variation" of Kunen's original argument.  All of these topics are covered extensively in \cite{MattHandbook}; but we only need a special case.  To describe the special case, we recall (see Kanamori~\cite{MR1994835}) that a sequence $\vec{U} = \langle U_\gamma \ | \ \gamma < \delta \rangle$ is called a \emph{$(\kappa,\delta)$ tower} iff each $U_\gamma$ is a normal measure on $P_\kappa(\gamma)$, and if $\gamma < \gamma'$ then
\[
U_{\gamma'} = \{ x \in P_\kappa(V_{\gamma'}) \ | \ x \cap V_\gamma \in U_\gamma  \}.
\]
Given such a tower, if $\gamma < \gamma'$ there is an elementary embedding $k_{\gamma,\gamma'}: \text{ult}(V,U_\gamma) \to \text{ult}(V,U_{\gamma'})$ such that $j_{\gamma'} = k_{\gamma,\gamma'} \circ j_\gamma$, where $j_{\gamma'}$ and $j_{\gamma}$ are the ultrapower maps for $U_{\gamma'}$ and $U_{\gamma}$, respectively.  This yields a directed system of ultrapower maps and a direct limit ultrapower embedding $j_{\vec{U}}: V \to_{\vec{U}} N_{\vec{U}}$ with critical point $\kappa$ such that $N_{\vec{U}}$ is wellfounded and in fact
\begin{equation}\label{eq_ClosedLessDelta}
N_{\vec{U}} \text{ is closed under } <\delta \text{ sequences}
\end{equation}
and in general $j_{\vec{U}}(\kappa) \ge \delta$.  If $j_{\vec{U}}(\kappa) = \delta$ then $\vec{U}$ is called an \emph{almost huge} $(\kappa,\delta)$ tower; in that case (by \eqref{eq_ClosedLessDelta}) notice that $N_{\vec{U}}$ is closed under $< j_{\vec{U}}\big(\text{crit}(j_{\vec{U}})\big)$ sequences.  We use the following fact:
\begin{fact}\label{fact_DerivedAlmostHuge}
Suppose $j:V \to N$ is a huge embedding with critical point $\kappa$; i.e.\ $j$ is elementary and $N$ is closed under $\delta:= j(\kappa)$ sequences.  Then there exists an almost huge $(\kappa,\delta)$-tower $\vec{U}$ and an elementary embedding $k: N_{\vec{U}} \to N$ such that $j = k \circ j_{\vec{U}}$ and $\text{crit}(k) = (\delta^+)^{N_{\vec{U}}}$.  Note in particular this implies that $k$ fixes $\delta$ and that $(\delta^+)^{N_{\vec{U}}}  < (\delta^+)^{N}$.
\end{fact}

An ideal $\mathcal{I}$ on a regular uncountable cardinal $\kappa$ is called \emph{normal} if it is closed under diagonal unions, and \emph{saturated} iff $\wp(\kappa)/\mathcal{I}$ has the $\kappa^+$-chain condition.  The following theorem is a special instance of the Duality Theorem of Foreman~\cite{MattHandbook} (see also Foreman-Komjath~\cite{MR2151585} or section 5.3 of Cox-Zeman~\cite{Cox_MALP} for more detailed proofs).   
\begin{theorem}\label{thm_Duality_AH}
Suppose $\vec{U}$ is an almost huge $(\kappa,\delta)$ tower and $j_{\vec{U}}: V \to_{\vec{U}} N_{\vec{U}}$ is the corresponding embedding.  Suppose $\mathbb{P} \subset V_\kappa$ is a $\kappa$-cc poset, $\dot{\mathbb{C}}_{\text{Silv}}(\kappa, < \delta)$ is the $\mathbb{P}$-name for the Silver collapse which turns $\delta$ into $\kappa^+$ (described in Section \ref{sec_Examples}), and in $N_{\vec{U}}$ there exists a regular embedding
\[
\iota: \mathbb{P}*\dot{\mathbb{C}}_{\text{Silv}}(\kappa, < \delta) \to j_{\vec{U}}(\mathbb{P})
\]
such that $\iota$ is the identity on $\mathbb{P}$.\footnote{I.e.\ $\iota(p,\dot{1}) = p$ for every $p \in \mathbb{P}$.}  Let $G*H$ be generic for $\mathbb{P}*\dot{\mathbb{C}}_{\text{Silv}}(\kappa, < \delta)$ over $V$.  Then in $V[G*H]$ there is a normal ideal $\mathcal{I}(j_{\vec{U}})$ on $\kappa$ such that $\wp(\kappa)/\mathcal{I}(j_{\vec{U}})$ is forcing equivalent to $j_{\vec{U}}(\mathbb{P})/\iota[G*H]$. 
\end{theorem}

We use one more fact about normal ideals:
\begin{fact}\label{fact_SaturatedImpliesCompleteBA}[\cite{MattHandbook}]
If $\mathcal{I}$ is a normal, saturated ideal on $\kappa$ then $\wp(\kappa)/\mathcal{I}$ is a complete boolean algebra. 
\end{fact}

\section{Layered posets}\label{sec_LayeredPosets}

\subsection{Basic theory of Layered posets}

A layered poset is, roughly, a poset which has many small regular suborders, where ``many" is typically taken to mean club or stationarily many.  There are two common definitions of stationarity for subsets of $P_\kappa(H)$, where $H$ is any set of cardinality at least $\kappa$:
\begin{itemize}
 \item Strong stationarity, as Definition 8.21 of Jech~\cite{MR1940513}.  A set $S \subseteq P_\kappa(H)$ is \emph{strongly stationary in $P_\kappa(H)$} iff $S \cap C \ne \emptyset$ for all $C \subseteq P_\kappa(H)$ which are cofinal in $(P_\kappa(H), \subset)$ and closed under increasing unions of length $<\kappa$.
 \item Shelah's Weak stationarity or generalized stationarity, e.g.\ as in Larson~\cite{MR2069032} and Foreman~\cite{MattHandbook}.  A set $S \subseteq P_\kappa(H)$ is \emph{weakly stationary in $P_\kappa(H)$} iff whenever $F: [H]^{<\omega} \to H$ then there is some $X \in S$ which is closed under $F$.  The set of all $X \subseteq H$ which are closed under $F$ is denoted by $C_F$.  
\end{itemize}
For $\kappa = \omega_1$ the notions ``strongly stationary in $P_\kappa(H)$" and ``weakly stationary in $P_\kappa(H)$" are equivalent, but for $\kappa \ge \omega_2$ they may differ, due to the presence or absence of Changs Conjecture; this is discussed in detail in Foreman~\cite{MattHandbook}.  
\begin{convention}\label{conv_JechStat}
Unless stated otherwise, in this paper ``stationary" will refer to Jech's Strong stationarity.\footnote{We could alternatively stick with weak stationarity, but require that our stationary subsets of $P_\kappa(H)$ concentrate on those $M \in P_\kappa(H)$ such that $M \cap \kappa \in \kappa$.}
\end{convention}
For example, in the proof of Lemma \ref{lem_LayeredImpliesKnaster}, it is essential that the stationary layering occurs on models $M$ such that $M \cap \kappa \in \kappa$; this happens automatically in any strongly stationary set, but if $\kappa > \omega_1$ it might fail for a weakly stationary set.

The following is the main definition:
\begin{definition}\label{def_Layer}
Let $\mathbb{P}$ be a poset, $\kappa$ a regular uncountable cardinal with $\kappa \le |\mathbb{P}|$, and
\begin{equation*}
\text{Reg}_{\kappa}(\mathbb{P}):= \{  X \in P_\kappa(\mathbb{P}) \ | \  X \text{ is a regular suborder of } \mathbb{P} \}.
\end{equation*}

\begin{enumerate}
  \item $\mathbb{P}$ is \emph{$\kappa$ stationarily layered} iff $\text{Reg}_{\kappa}(\mathbb{P})$ is stationary in $P_\kappa(\mathbb{P})$.
  \item $\mathbb{P}$ is \emph{$\kappa$ club layered} iff $\text{Reg}_{\kappa}(\mathbb{P})$ contains a club (i.e.\ $P_\kappa(\mathbb{P}) - \text{Reg}_{\kappa}(\mathbb{P})$ is not stationary).
  \item $\mathbb{P}$ is \emph{$\kappa$ cofinally layered} iff $\text{Reg}_{\kappa}(\mathbb{P})$ is cofinal in $\big( P_\kappa(\mathbb{P}), \subseteq \big)$.
\end{enumerate}

Recall Convention \ref{conv_JechStat}; we are referring to Jech's strong stationarity here.
\end{definition}

\begin{remark}\label{rem_TFAE_statlayer}
By standard lifting and projecting of stationary sets (see Theorem 8.27 of \cite{MR1940513}), the following are equivalent:
\begin{enumerate}
 \item $\mathbb{P}$ is $\kappa$ stationarily layered;
 \item The set
 \begin{equation*}
 \{ X \in P_\kappa(H) \ | \  X \cap \mathbb{P} \text{ is a regular suborder of } \mathbb{P} \}
 \end{equation*}
 is stationary in $P_\kappa(H)$, where $H$ is any set such that $\mathbb{P} \subset H$. 
\end{enumerate}

The analogous equivalence also holds if we replace ``stationary" by ``contains a club" (and if we replace ``stationary" by ``weakly stationary").
\end{remark}

\begin{definition}
Assume $\mathbb{P} \subset H$ and $S \subseteq P_\kappa(H)$ is stationary.  We say that $\mathbb{P}$ is layered a.e.\ on $S$ iff $\{ M \cap \mathbb{P} \ | \ M \in S \} - \text{Reg}_{\kappa}(\mathbb{P})$ is nonstationary.  Equivalently, $M \cap \mathbb{P}$ is regular in $\mathbb{P}$ for all but nonstationarily many $M \in S$.
\end{definition}

A special case of the following lemma that appears in the literature is Shelah's work on \emph{layered ideals}.  
\begin{lemma}\label{lem_StatLayerImpliesKappaCC}[Shelah (see \cite{MattHandbook})]
If $\mathbb{P}$ is $\kappa$ weakly stationarily layered---i.e.\ if the set of regular suborders of $\mathbb{P}$ is weakly stationary in $P_\kappa(\mathbb{P})$---then $\mathbb{P}$ is $\kappa$-cc.
\end{lemma}
\begin{proof}
Let $A$ be a maximal antichain of $\mathbb{P}$.  Let $\theta$ be a sufficiently large regular cardinal so that $A, \mathbb{P} \in H_\theta$.  By Remark \ref{rem_TFAE_statlayer} (the ``weakly stationary" version), there is an $X \prec (H_\theta , \in, \{ \mathbb{P}, A \})$ (note that $A \in X$) with $|X|<\kappa$ such that $X \cap \mathbb{P}$ is a regular suborder of $\mathbb{P}$.  Because $X \prec (H_\theta , \in, \{ \mathbb{P}, A \})$, then $A \cap X$ is a maximal antichain in $X \cap \mathbb{P}$.  Because $X \cap \mathbb{P}$ is regular in $\mathbb{P}$, then $A \cap X$ is a maximal antichain in $\mathbb{P}$.  It follows that $A \cap X = A$; i.e.\ that $A$ is a subset of the $<\kappa$-sized set $X$.  So $|A| < \kappa$.
\end{proof}

For posets of size $\kappa$, it is often more natural to work with \emph{filtrations} of the poset to determine if the poset is layered.  A filtration of a set $X$ is a $\subset$-increasing, $\subset$-continuous sequence $\langle X_\xi \ | \ \xi < \text{cf}(|X|) \rangle$ such that $|X_\xi| < |X|$ for all $\xi < |X|$, and $X = \bigcup_{\xi < \text{cf}(|X|)} X_\xi$.  If $\text{cf}(|X|)$ is uncountable then any two filtrations agree on a club subset of $\text{cf}(|X|)$.  The proof of the following is standard so we omit it.
\begin{lemma}\label{lem_FiltrationLayering}
Suppose $\mathbb{P}$ is a poset of size $\kappa$, where $\kappa$ is regular and uncountable.  The following are equivalent:
\begin{enumerate}
 \item $\mathbb{P}$ is $\kappa$ stationarily layered
 \item For some filtration $\langle \mathbb{Q}_\alpha \ |  \ \alpha < \kappa \rangle$ of $\mathbb{P}$, there are stationarily many $\alpha < \kappa$ such that $\mathbb{Q}_\alpha$ is a regular suborder of $\mathbb{P}$.
 \item For every filtration $\langle \mathbb{Q}_\alpha \ |  \ \alpha < \kappa \rangle$ of $\mathbb{P}$, there are stationarily many $\alpha < \kappa$ such that $\mathbb{Q}_\alpha$ is a regular suborder of $\mathbb{P}$.
\end{enumerate} 
\end{lemma}

If $M \prec (H_\theta,\in,\mathbb{P})$ and $p \in \mathbb{P}$, then $p$ is called an $(M,\mathbb{P})$-master condition iff $p \Vdash \check{M}[\dot{G}] \cap \text{ORD} = \check{M} \cap \text{ORD}$.  The following correspondence between $\kappa$-cc and master conditions is due to Mekler:
\begin{lemma}[Mekler~\cite{MR758934}]\label{lem_Mekler}
Given a poset $\mathbb{P}$ and a regular uncountable cardinal $\kappa$, the following are equivalent:
\begin{enumerate}
 \item $\mathbb{P}$ is $\kappa$-cc;
 \item There are stationarily many $M  \in P_\kappa(H_\theta)$ such that $1_{\mathbb{P}}$ is an $(M, \mathbb{P})$-master condition.
\end{enumerate}

If $\kappa = \omega_1$ then these are also equivalent to:  there are club many $M  \in P_\kappa(H_\theta)$ such that $1_{\mathbb{P}}$ is an $(M, \mathbb{P})$-master condition.
\end{lemma}

Just as Lemma \ref{lem_Mekler} characterizes $\kappa$-cc  in terms of master conditions, $\kappa$-stationary layering can be expressed in terms of \emph{strong} master conditions, as defined by Mitchell~\cite{MR2279659}.  If $M$ is a model, $\mathbb{P}$ is a poset, and $p \in \mathbb{P}$, Mitchell~\cite{MR2452816} defined $p$ to be an $(M,\mathbb{P})$ strong master condition iff every $p' \le p$ has a reduction in $M \cap \mathbb{P}$.  Note that if $\mathbb{P} \subset H_\theta$ and $X \prec (H_\theta, \in, \mathbb{P})$, then by elementarity $X \cap \mathbb{P}$ is a suborder of $\mathbb{P}$ (i.e.\ we have order and $\perp$ preservation in both directions).  The following observation is straightforward, using the ``reduct" characterization of regularity in Lemma \ref{lem_ReductCharactOfReg}:
\begin{observation}\label{obs_RegSub_StrongMaster}
If $X \prec (H_\theta, \in \mathbb{P})$, the following are equivalent:
\begin{itemize}
 \item $X \cap \mathbb{P}$ is a regular suborder of $\mathbb{P}$.
 \item $1_{\mathbb{P}}$ is an $(X,\mathbb{P})$-strong master condition (in the sense of Mitchell~\cite{MR2452816}).
\end{itemize}
\end{observation}

\begin{corollary}\label{cor_TFAE_StatLayer_1StrongMaster}
Given a poset $\mathbb{P}$ and a regular uncountable cardinal $\kappa$, the following are equivalent:
\begin{enumerate}
 \item $\mathbb{P}$ is $\kappa$ stationarily layered;
 \item There are stationarily many $M  \in P_\kappa(H_\theta)$ such that $1_{\mathbb{P}}$ is a $(M, \mathbb{P})$ \emph{strong} master condition.
\end{enumerate}
\end{corollary}
\begin{proof}
This follows immediately Observation \ref{obs_RegSub_StrongMaster}.
\end{proof}

\subsection{Examples of layered posets}\label{sec_Examples}

Many commonly used posets are stationarily  layered under mild cardinal arithmetic assumptions; in fact they are often layered a.e.\ on some natural stationary set $S$ (i.e.\ ``club-layered relative to $S$").  Here are a few common examples of $\kappa$-stationarily layered posets, all under the assumption that $\kappa$ is (strongly) inaccessible.  Much less than inaccessibility is needed for most of these, given the right cardinal arithmetic assumptions. 
\begin{enumerate}
 \item Trivially, any poset of density $<\kappa$ is $\kappa$-club layered;
 \item Any poset as in the hypothesis of Corollary \ref{cor_IterateSmallPosets} is $\kappa$ stationarily layered. 
 \item If $\kappa$ is inaccessible then the Levy collapse $\text{Col}(\mu, < \kappa)$ is layered a.e.\ on the stationary set of $<\mu$-closed models in $P_\kappa(H_{\kappa})$ (recall we are using strong stationarity here, so almost every such model has transitive intersection with $\kappa$).  If $\mu = \omega$ then it is club-layered on $P_{\omega_1}(V_\kappa)$.
 \item If $\kappa$ is inaccessible then the Silver Collapse $\mathbb{C}_{\text{Silv}}(\mu, < \kappa)$ that turns $\kappa$ into $\mu^+$ is layered a.e.\ on the stationary set of $\mu$-closed elements of $P_\kappa(H_{\kappa})$.  The poset $\mathbb{C}_{\text{Silv}}(\mu, < \kappa)$ is the $\mu$-support product of $\{ \text{Col}(\mu,\eta) \ | \ \eta < \kappa \}$ with bounded domain.  That is, conditions are partial functions $p: \mu \times \kappa \to \kappa$ such that $|p| \le \mu$, $p(\alpha,\eta) < \eta$ whenever $(\alpha,\eta) \in \text{dom}(p)$, and there is some $\xi_p < \mu$ such that whenever $(\alpha,\eta) \in \text{dom}(p)$ then $\alpha < \xi$.
 \item\label{item_EastonCollapse} If $\kappa$ is Mahlo then the Easton Collapse $\mathbb{E}(\mu,<\kappa)$ is layered a.e.\ on the stationary set of $W \in P_\kappa(V_\kappa)$ such that $|W| = W \cap \kappa \in \kappa$, $\mu \in W$, and ${}^{<|W|} W \subset W$.  Conditions are partial functions $p: \mu \times \kappa \to \kappa$ such that $p(\alpha,\eta) < \eta$ whenever $(\alpha,\eta) \in \text{dom}(p)$, there is some $\xi < \mu$ such that $\alpha < \xi$ whenever $(\alpha, \eta) \in \text{dom}(p)$, and the support of $p$---i.e.\ the set $\{ \eta < \kappa \ | \ \exists \alpha \ (\alpha,\eta) \in \text{dom}(p) \}$---is \emph{Easton above $\mu$}.  A set $X \subset \kappa$ is Easton above $\mu$ iff $|X \cap \gamma|<\gamma$ whenever $\gamma \in (\mu,\kappa)$ is regular. 
 \item The poset $\mathbb{B}(\mu,\alpha, \kappa)$ from Foreman-Komjath~\cite{MR2151585} is layered a.e.\ on the stationary set of $<\alpha$-closed elements of $P_\kappa(V_\kappa)$.  Here $\mu < \alpha < \kappa$ are regular cardinals and conditions are elements of the $<\alpha$-support product of $\kappa$-many copies of $\text{Col}(\mu,\alpha)$ which have a bounded domain.  More precisely, let $\mathbb{C}_i:= \text{Col}(\mu,\alpha)$ for each $i < \kappa$ (so these are $\kappa$-many copies of $\text{Col}(\mu,\alpha)$).  A condition is a $< \kappa$-supported function $p \in \prod_{i<\kappa} \mathbb{C}_i$ such that there is some $\xi_p < \mu$ such that $\text{dom}(p(i)) \subseteq \xi_p$ for all $i \in \text{supp}(p)$.  
 \item Any $\kappa$-cc poset---and any $\mu$-support product of $\kappa$-cc posets where $\mu < \kappa$---is $\kappa$ stationarily layered, if $\kappa$ is weakly compact, by \cite{Cox_Luecke}.
\end{enumerate} 

We prove that the Easton collapse $\mathbb{E}(\mu, < \kappa)$ is layered on the stationary set indicated above; the proofs that the other posets are layered is similar (with the exception of the weakly compact example, which is proved in \cite{Cox_Luecke}).  Pick any $W \prec (H_{\kappa^+}, \in)$ such that $\gamma_W:=|W|=W \cap \kappa \in \kappa$ and ${}^{<|W|} W \subset W$.  Let $p \in \mathbb{E}(\mu,<\kappa)$, and let $S$ be the support of $p$, which is Easton above $\mu$ (as defined in \ref{item_EastonCollapse} above).  Because $S$ is Easton, $\gamma_W$ is regular (inaccessible), and $W$ is closed under $<\gamma_W$ sequences, then $p \restriction W$ is an element of $W$.  Note that because $W \cap \kappa \in \kappa$ then $p \restriction W$ is the same as $p \restriction  (\mu \times \gamma_W)$.  We need to prove that it is a reduct.  Suppose $p' \in W$ is a condition extending $p \restriction W$.  Then clearly $p' \cup p$ is a function, because $p'$, being in $W$, does not make any commitments for those $\eta \notin W$.  Let $\xi_{p'}$ be the uniform bound for $p'$'s domains, and $\xi_{p}$ be the uniform bound for $p$'s domains.  Then $\text{max}(\xi_{p'}, \xi_{p})$ is a uniform bound for $p \cup p'$.  So $p \cup p'$ is a condition.  

\begin{remark}
Naturally-defined layered posets as in the examples listed above typically have $\kappa$-stationarily layered products, simply because the stationary sets witnessing layering have stationary overlap.  For example, Foreman-Komjath~\cite{MR2151585} consider products of the form
\begin{equation*}
\mathbb{S}(\alpha, < \kappa) \times \mathbb{B}(\mu, \alpha, \kappa)
\end{equation*}
where $\kappa$ is inaccessible.  As remarked above, $\mathbb{S}(\alpha, < \kappa)$ is layered a.e.\ on the set of $\alpha$-closed models, and $\mathbb{B}(\mu, \alpha, \kappa)$ is layered a.e.\ on the set of $<\alpha$-closed models.  Because the intersection of these stationary sets is the set of $\alpha$-closed models, which is stationary, it follows that their product is layered a.e.\ on the set of $\alpha$-closed models; this and other issues surrounding products are addressed in detail in Cox-L\"ucke~\cite{Cox_Luecke}.  It also follows from the more general fact that if $p_i$ is a strong master condition for $(M, \mathbb{P}_i)$ for $i \in \{ 0, 1 \}$ then $(p_0, p_1)$ is a strong master condition for $(M,  \mathbb{P}_0 \times \mathbb{P}_1)$; see Mitchell~\cite{MR2452816}.
\end{remark}

\subsection{Relation to Knaster property}\label{sec_RelationToKnaster}

A poset is called $\kappa$-Knaster if every $\kappa$-sized collection of conditions can be refined to a $\kappa$-sized collection of pairwise compatible conditions.  Clearly $\kappa$-Knaster implies $\kappa$-cc.  The main feature of the Knaster property used in the literature is:

\begin{fact}\label{fact_KnasterProduct}[see \cite{MR2768691}]
If $\mathbb{P}$ is $\kappa$-Knaster and $\mathbb{Q}$ is $\kappa$-cc, then $\mathbb{P} \times \mathbb{Q}$ is $\kappa$-cc.
\end{fact}

We now prove Lemma \ref{lem_LayeredImpliesKnaster}, that $\kappa$-stationary layering implies $\kappa$-Knaster.  The proof is modeled after the proof of Theorem 4.1 of Soukup~\cite{MR2800978}.

Let $\mathbb{P}$ be $\kappa$ stationarily layered, and let $X$ be a $\kappa$-sized collection of conditions in $\mathbb{P}$.  By the layering assumption of $\mathbb{P}$, there is an $M \prec (H_\theta, \in, \mathbb{P})$ such that $X \in M$, $|M|<\kappa$, $M \cap \kappa\in \kappa$, and $M \cap \mathbb{P}$ is a regular suborder of $\mathbb{P}$.  Because $|M|<\kappa=|X|$, there is some $p \in X - M$.  Let $p|M$ be any reduct of $p$ into $M \cap \mathbb{P}$.  Let $Z$ be a maximal subset of $X$ with the property that $Z$ consists of pairwise compatible conditions, and every member of $Z$ is compatible with $p|M$.  Since $X$ and $p|M$ are members of $M$, we can without loss of generality take $Z \in M$.  We claim that $|Z| = \kappa$.  If not then because $Z \in M$ and $|Z|<\kappa$, it follows that $|Z| \in M \cap \kappa$.  And because $M \cap \kappa$ is transitive,\footnote{Note the essential use of transitivity of $M \cap \kappa$ in this part of the argument.}
\begin{equation}\label{eq_Z_subset_M}
Z \subset M.
\end{equation}
It follows that $Z \subsetneq Z \cup \{ p \}$.  Now $p$ is clearly compatible with $p|M$, and $Z \cup \{ p \}$ is a subset of $X$; thus we will obtain a contradiction to maximality of $Z$ if we can prove that $Z \cup \{ p \}$ is pairwise compatible.  Let $r \in Z$; then \eqref{eq_Z_subset_M} implies $r \in M \cap \mathbb{P}$ and, by the requirements on $Z$, we have $r \parallel_{\mathbb{P}} p|M$.  Since $r$ and $p|M$ are in $M$ then this compatibility is witnessed by some $r' \le r, p|M$ such that $r' \in M$.  Because $r' \le p|M$ it follows that $r' \parallel_{\mathbb{P}} p$, which implies that $r \parallel_{\mathbb{P}} p$.  This proves that $Z \cup \{ p \}$ is pairwise compatible, completing the proof.

\begin{remark}
It is natural to wonder if the implication of Lemma \ref{lem_LayeredImpliesKnaster} can be reversed.  If $\kappa$ is weakly compact then the answer is yes; and for non-weakly compact cardinals the answer is, in general, no; see Cox-L\"ucke~\cite{Cox_Luecke}.
\end{remark}

\subsection{Layering and elementary embeddings}\label{sec_LayeringElementaryEmbeddings}

The following simple but useful lemma will be used in the proof of Theorem \ref{thm_WeakCompact_Kunen_MainThm_Imprecise} (and also in \cite{Cox_Luecke}):
\begin{lemma}\label{lem_EmbeddingsAndKappaCC}
Suppose 
\begin{itemize}
 \item $j: H \to N$ is a $\Sigma_0$ elementary embedding of $ZFC^-$ models
 \item $\kappa= \text{crit}(j)$
 \item $j[H] \in N$ (equivalently $j \in N$)
 \item $\mathbb{Q} \in H$ is a poset
 \item $N \models \mathbb{Q}$ is $\kappa$-cc
 \item $([\mathbb{Q}]^{<\kappa})^H = ([\mathbb{Q}]^{<\kappa})^N$.
\end{itemize}
Then $j[H] \cap j(\mathbb{Q})$ is a regular suborder of $j(\mathbb{Q})$. 
\end{lemma}
\begin{proof}
Because $j[H] \in N$ and ``$X$ is a regular suborder of $Y$" is a $\Sigma_0$ property (using the reduct characterization), it suffices to prove that $N$ believes $j[H] \cap j(\mathbb{Q})$ is regular in $j(\mathbb{Q})$.  Clearly $j \restriction \mathbb{Q}: \mathbb{Q} \to j(\mathbb{Q})$ is order and incompatibility preserving, so we only need to check that whenever $A \in N$ is a maximal antichain in $\mathbb{Q}$ then $j[A]$ is maximal in $j(\mathbb{Q})$.  So let $A \in N$ be maximal antichain in $\mathbb{Q}$.  Because $([\mathbb{Q}]^{<\kappa})^H = ([\mathbb{Q}]^{<\kappa})^N$ and $N \models$ $\mathbb{Q}$ is $\kappa$-cc" then $A \in H$ and has size $<\kappa$ in $H$.  Since $\text{crit}(j) = \kappa$ this implies that $j(A) = j[A]$.  But by elementarity, $N$ believes that $j(A)$ is maximal in $j(\mathbb{Q})$, completing the proof. 
\end{proof}

\section{Coherent, conservative systems of reduction operations}\label{sec_InductReductions}

The iteration we describe in the main theorem isn't presented in the usual way, so we will make the following ad-hoc definition; it is really a definition about how a poset is presented:
\begin{definition}\label{def_GeneralizedIteration}
A sequence of posets $\langle \mathbb{L}_\alpha \ | \ \alpha < \theta \rangle$ will be called a \emph{generalized iteration} iff:
\begin{itemize}
 \item Each $\mathbb{L}_\alpha$ consists of partial functions on $\alpha$;
 \item The restriction map $\pi^\beta_\alpha$ is a forcing projection from $\mathbb{L}_\beta \to \mathbb{L}_\alpha$ for all $\beta \ge \alpha$.  That is, $\pi^\beta_\alpha$ is order preserving, and whenever $q \le \pi^\beta_\alpha(p)$ then there is some $p' \le p$ such that $\pi^\beta_\alpha(p') \le q$.
 \item The identity map is a regular embedding from $\mathbb{L}_\alpha \to \mathbb{L}_\beta$ for all $\alpha \le \beta$.
\end{itemize}
\end{definition}

The following definition very important for the main theorem.  If $\mathbb{L}$ is a suborder of $\mathbb{R}$, a map $\pi: \mathbb{R} \to \mathbb{L}$ is called a \emph{reduction operation} iff $\pi(r)$ is a reduct of $r$ into $\mathbb{L}$ for all $r \in \mathbb{R}$;\footnote{See Section \ref{sec_Prelims} for a discussion of reducts.} i.e.\ whenever $\ell \le \pi(r)$ and $\ell \in \mathbb{L}$ then $\ell$ is compatible with $r$.\footnote{Note that this is \textbf{not} the same as a forcing projection; a forcing projection would require instead that there is some $r' \le r$ such that $\pi(r') \le \ell$.}
\begin{definition}\label{def_ConservCoherent}
Assume $\vec{\mathbb{L}} = \langle \mathbb{L}_\alpha \ | \ \alpha < \theta \rangle$ and $\vec{\mathbb{R}} = \langle \mathbb{L}_\alpha \ | \ \alpha < \theta \rangle$ are two generalized iterations of the same length $\theta$, and assume that $\mathbb{L}_\alpha$ is a regular suborder of $\mathbb{R}_\alpha$ for all $\alpha < \theta$.  A system of reduction maps $\langle \pi_\alpha: \mathbb{R}_\alpha \to \mathbb{L}_\alpha \ | \ \alpha < \theta \rangle$ will be called:
\begin{itemize}
 \item \emph{coherent} iff $\pi_\alpha(r \restriction \alpha) = \pi_\beta(r) \restriction \alpha$ for all $\alpha \le \beta < \theta$ and all $r \in \mathbb{R}_\beta$;
 \item \emph{conservative} iff whenever:
 \begin{itemize}
  \item $\alpha \le \beta$;
  \item $\ell_\beta \le \pi_\beta(r_\beta)$
  \item $r'_\alpha$ is some condition in $\mathbb{R}_\alpha$ witnessing compatibility of $\ell_\beta \restriction \alpha$ with $r_\beta \restriction \alpha$;  
 \end{itemize}
 then there is some $r'_\beta \in \mathbb{R}_\beta$ such that $r'_\beta \restriction \alpha = r'_\alpha$, and $r'_\beta$ witnesses compatibility of $\ell_\beta$ with $r_\beta$.
\end{itemize}
\end{definition}

The conservativity requirement is used to pass inverse limit stages of inductively defined reduction mappings.  Roughly, it says that if you're checking compatibility (in a proof that some condition is a reduct), you can do it without strengthening what you did earlier.

\begin{lemma}\label{lem_WhatToDoAtLimits}
Assume $\vec{\mathbb{L}}$ and $\vec{\mathbb{R}}$ are generalized iterations of limit length $\theta$, and that $\vec{\pi} = \langle \pi_\alpha: \mathbb{R}_\alpha \to \mathbb{L}_\alpha \ | \ \alpha < \theta \rangle$ is a coherent, conservative system of reduction maps.  Assume $\vec{\mathbb{L}}$ and $\vec{\mathbb{R}}$ use some mix of inverse and direct limits, and that their limit scheme is the same.  More precisely, for each limit ordinal $\alpha$, either 
\begin{enumerate*}
 \item $\mathbb{L}_\alpha$ and $\mathbb{R}_\alpha$ are both direct limits; or
 \item $\mathbb{L}_\alpha$ and $\mathbb{R}_\alpha$ are both inverse limits. 
\end{enumerate*}

Let $\mathbb{L}_\theta$ and $\mathbb{R}_\theta$ be limits of $\vec{\mathbb{L}}$ and $\vec{\mathbb{R}}$, respectively, where either both are direct limits or both are inverse limits.  Then there is a reduction map $\pi_\theta: \mathbb{R}_\theta \to \mathbb{L}_\theta$ such that $\vec{\pi} \cup \{( \theta, \pi_\theta ) \}$ is a coherent, conservative system of reduction maps.
\end{lemma}
\begin{proof}
We prove the lemma assuming that both $\mathbb{L}_\theta$ and $\mathbb{R}_\theta$ are inverse limits, as this is the hardest case.

Notice that the coherency of $\vec{\pi}$ ensures that:
\begin{equation}\label{eq_WhatLimitsMustLookLike}
\forall \gamma \in \theta \cap \text{Lim} \ \  \forall r_\gamma \in \mathbb{R}_\gamma \ \ \pi_\gamma(r_\gamma) = \bigcup_{\alpha < \gamma} \pi_\alpha(r_\gamma \restriction \alpha).
\end{equation}
Also, the coherency requirement that we desire of $\pi_\theta$ compels us to define 
\begin{equation*}
\pi_\theta(r):=  \bigcup_{\alpha < \theta} \pi_\alpha(r \restriction \alpha).
\end{equation*}
That $\pi_\theta(r)$ is a condition in $\mathbb{L}_\theta$ follows from the facts that:
\begin{itemize}
 \item Coherency of $\vec{\pi}\restriction \theta$ ensures that $\pi_\beta(r \restriction \beta)$ end-extends $\pi_\alpha(r \restriction \alpha)$ whenever $\alpha \le \beta$
 \item $\mathbb{L}_\theta$ and $\mathbb{R}_\theta$ use the same limit type.
\end{itemize} 

We must prove that conservativity is preserved (the proof will also show that $\pi_\theta$ is a reduction mapping).  So assume $\alpha < \theta$ and
\begin{itemize}
  \item $\ell_\theta \le \pi_\theta(r_\theta)$
  \item $r'_\alpha$ is some condition in $\mathbb{R}_\alpha$ witnessing compatibility of $\ell_\theta \restriction \alpha$ with $r_\theta \restriction \alpha$;  
 \end{itemize}
We need to find some $r'_\theta \in \mathbb{R}_\theta$ such that $r'_\theta \restriction \alpha = r'_\alpha$, and $r'_\theta$ witnesses compatibility of $\ell_\theta$ with $r_\theta$.  We construct such an $r'_\theta$ using the assumption that $\vec{\pi}$ is conservative, together with \eqref{eq_WhatLimitsMustLookLike}.  Namely, recursively define a sequence $\langle r'_\beta \ | \ \alpha \le \beta \le \theta \rangle$ as follows:
\begin{itemize}
 \item If $r'_\beta$ is defined, let $r'_{\beta+1}$ be some condition in $\mathbb{R}_{\beta+1}$ given by the conservativity between $\beta$ and $\beta + 1$.  That is:
 \begin{itemize}
  \item  $r'_{\beta+1} \restriction \beta = r'_\beta$ 
  \item $r'_{\beta+1}$ witnesses compatibility of $\ell_\theta \restriction (\beta+1)$ with $r_\theta \restriction (\beta+1)$.  If index $\theta$ is a trivial coordinate for both $\ell_\theta$ and $r_\theta$---i.e.\ if  $\ell_\theta \restriction (\beta+1) = \ell_\theta \restriction \beta$ and $r_\theta \restriction (\beta+1) = r_\theta \restriction \beta$---then we require that the $\theta$-th coordinate of $r'_{\beta+1}$ is also trivial; more precisely that $r'_{\beta+1} = r'_\beta$.   
 \end{itemize}
 
 \item If $\gamma$ is limit, let $r'_\gamma = \bigcup_{\beta \in [\alpha, \gamma)} r'_\beta$.
\end{itemize}

The triviality requirement at successor steps of the construction, together with the assumption that $\vec{\mathbb{L}}$ and $\vec{\mathbb{R}}$ use the same limit scheme, ensures that each $r'_\beta$ is a condition in $\mathbb{R}_\beta$.  It clear that $r'_\beta = r'_{\beta'} \restriction \beta$ whenever $\beta \le \beta'$, and that $r'_\theta$ is as required. 
\end{proof}

Though Lemma \ref{lem_WhatToDoAtLimits} assumes that the iterations use some mix of inverse and direct limits, we suspect a similar fact holds for many common iteration schemes (e.g.\ RCS iterations).

\section{Universal Kunen iterations}\label{sec_KunenUniversal}

The universal Kunen iteration was introduced by Kunen~\cite{MR495118} in the first consistency proof of a saturated ideal on $\omega_1$, and is a useful tool when huge embeddings are involved.  Roughly, one has a large cardinal $\kappa$ and defines some sort of iteration of length $\kappa$, where successor steps are handled as follows:  given that $0< \alpha < \kappa$ and $\mathbb{P}_\alpha$ has been defined, if $V_\alpha \cap \mathbb{P}_\alpha$ is a regular suborder of $\mathbb{P}_\alpha$, set $\mathbb{P}_{\alpha+1} = \mathbb{P}_\alpha * \dot{\mathbb{Q}}_\alpha$ where $\dot{\mathbb{Q}}_\alpha$ is some poset \textbf{from the inner model} $V[\dot{G}_\alpha \cap V_\alpha]$; otherwise set $\mathbb{P}_{\alpha+1} = \mathbb{P}_\alpha * \{ 1 \}$.  In the applications below it will be helpful to describe the Kunen iteration a bit differently from typical iterations (but these will still be iterations in the general sense of iteration, that is a commutative system of forcing projections).  The scheme below was roughly described in Foreman~\cite{MattHandbook}.
\begin{definition}\label{def_KunenIteration}
Let $\kappa$ be an uncountable regular cardinal.  A sequence $\langle \mathbb{P}_\alpha, \dot{Q}_\beta \ | \ \alpha \le \kappa, \beta < \kappa \rangle$ will be called a \emph{universal Kunen iteration} iff:
\begin{enumerate}
 \item $\langle \mathbb{P}_\alpha \ | \ \alpha \le \kappa \rangle$ is a generalized iteration as in Definition \ref{def_GeneralizedIteration};

 \item For each $\alpha < \kappa$, the poset $\mathbb{P}_{\alpha+1}$ has the following form.  If $V_\alpha \cap \mathbb{P}_\alpha$ is \textbf{not} a regular suborder of $\mathbb{P}_\alpha$ then we will call $\alpha$ a \emph{passive stage} and set $\mathbb{P}_{\alpha+1} = \mathbb{P}_\alpha$.  Otherwise---i.e.\ if $V_\alpha \cap \mathbb{P}_\alpha$ is a regular suborder of $\mathbb{P}_\alpha$---then we will call $\alpha$ an \emph{active stage} and require that $\dot{Q}_\alpha$ is a $V_\alpha \cap \mathbb{P}_\alpha$-name for a poset,\footnote{We emphasize that $\dot{Q}_\alpha$ is a $V_\alpha \cap \mathbb{P}_\alpha$-name, not merely a $\mathbb{P}_\alpha$-name.} and $\mathbb{P}_{\alpha+1}$ consists of all partial functions $f$ on $\alpha+1$ with the following properties:
 \begin{itemize}
  \item $f \restriction \alpha \in \mathbb{P}_\alpha$;
  \item If $\alpha \in \text{dom}(f)$ then $f(\alpha) \in V^{V_\alpha \cap \mathbb{P}_\alpha}$ (notice this is a $V_\alpha \cap \mathbb{P}_\alpha$-name, not merely a $\mathbb{P}_\alpha$-name) and $\Vdash_{V_\alpha \cap \mathbb{P}_\alpha} f(\alpha) \in \dot{Q}_\alpha$.
 \end{itemize}
where the partial ordering is given by: 
\begin{align*}
f_1 \le_{\mathbb{P}_{\alpha+1}} f_0 \ :\equiv \  f_1 \restriction \alpha \le_{\mathbb{P}_\alpha} f_0 \restriction \alpha \text{ and if } \alpha \in \text{dom}(f_0) \text{ then } \\
\alpha \in \text{dom}(f_1) \text{ and } f_1 \restriction \alpha \Vdash^{\text{Convention}}_{\mathbb{P}_\alpha} f_1(\alpha)  \le_{\dot{Q}_\alpha} f_0(\alpha)
\end{align*}
where we are using Convention \ref{conv_ViewNames}, together with the assumption that $V_\alpha \cap \mathbb{P}_\alpha$ is a regular suborder of $\mathbb{P}_\alpha$, to view the statement ``$f_1(\alpha)  \le_{\dot{Q}_\alpha} f_0(\alpha)$'' as a statement in the language of $\mathbb{P}_\alpha$.

\end{enumerate}
\end{definition}

As usual, we can arrange that $\mathbb{P}$ is a set (rather than a proper class) by choosing a single representative of the various names in the definition (modulo equivalence of names).  We will insist that such a representative is always chosen of minimal rank.  The following is a standard fact:
\begin{fact}\label{fact_NamesInModel}
Suppose $\gamma$ is regular and $\mathbb{R} \in H_\gamma$ is a poset.  If $\tau$ is a $\mathbb{R}$-name which is forced to be in $H_\gamma[\dot{G}_{\mathbb{R}}]$, then there is a $\mathbb{R}$-name $\tau' \in H_\gamma$ which is $\mathbb{R}$-equivalent to $\tau$ (i.e.\ $\Vdash_{\mathbb{R}} \tau = \tau'$).
\end{fact}

\begin{lemma}\label{lem_SubsetVkappa}
Suppose $\langle \mathbb{P}_\alpha, \dot{Q}_\beta \ | \ \alpha  \le \kappa, \beta < \kappa \rangle$ is a universal Kunen iteration as in Definition \ref{def_KunenIteration} which uses some mix of inverse and direct limits, and that for every active $\alpha < \kappa$:
\begin{equation}\label{eq_AssumeInVkappa}
\Vdash_{V_\alpha \cap \mathbb{P}_\alpha} \dot{Q}_\alpha \subseteq V_\kappa[\dot{g}_\alpha]
\end{equation}
where $\dot{g}_\alpha$ is the canonical $V_\alpha \cap \mathbb{P}_\alpha$-name for its generic object.  Also assume that $\mathbb{P}_\kappa$ is a direct limit.  Then:
\begin{enumerate}
 \item\label{item_CanWLOGassumeVkappa} $\mathbb{P}_\kappa \subset V_\kappa$.
 \item\label{item_IfYouHaveEmbedding} If $j: V \to N$ is an elementary embedding with critical point $\kappa$ then:
 \begin{itemize}
  \item $j[\mathbb{P}_\kappa] = \mathbb{P}_\kappa = j(\vec{\mathbb{P}})_\kappa \cap V_\kappa$;
  \item If $\mathbb{P}_\kappa$ is $\kappa$-cc then $\kappa$ is an active stage of $j(\vec{\mathbb{P}})$ (and thus by reflection there are many active stages of $\vec{\mathbb{P}}$).
 \end{itemize}
\end{enumerate}
\end{lemma}
\begin{proof}
For item \ref{item_CanWLOGassumeVkappa}:  because $\mathbb{P}_\kappa$ is a direct limit then it suffices to prove that each $\mathbb{P}_\alpha \subset V_\kappa$ (for $\alpha < \kappa$); note that given the setup of Definition \ref{def_KunenIteration}, a direct limit is simply a union.  

We prove by induction that $\mathbb{P}_\alpha \subseteq V_\kappa$ for every $\alpha < \kappa$.\footnote{Recall that we use the convention that all names for conditions are taken to be of minimal rank.}  For limit stage $\alpha$ this follows simply from the induction hypothesis and inaccessibility of $\kappa$ (the latter is needed if $\mathbb{P}_\alpha$ is, say, an inverse limit).  Now suppose $\mathbb{P}_\alpha \subseteq V_\kappa$ and let $f \in \mathbb{P}_{\alpha+1}$.  We need to prove that $f \in V_\kappa$; because $f \restriction \alpha \in \mathbb{P}_\alpha \subset V_\kappa$ by the induction assumption, it only remains to show that $f(\alpha) \in V_\kappa$.  This follows from assumption \eqref{eq_AssumeInVkappa} and Fact \ref{fact_NamesInModel} (together with our convention that a representative for a name is always chosen of minimal rank).  The proof that $V_\kappa \cap j(\vec{\mathbb{P}})_\kappa \subset \mathbb{P}_\kappa$ is proved similarly.

For item \ref{item_IfYouHaveEmbedding}: by item \ref{item_CanWLOGassumeVkappa} $j \restriction \mathbb{P}_\kappa$ is just the identity map on $\mathbb{P}_\kappa$. If $f \in \mathbb{P}_\kappa$ then because $\mathbb{P}_\kappa$ is a direct limit, $f \in \mathbb{P}_\alpha$ for some $\alpha < \kappa$.   So $j(f) = f$.  Pick any $\beta < \alpha$ in the domain of $f$.  Then $f \restriction \beta \in \mathbb{P}_\beta$ and $V_\beta \cap \mathbb{P}_\beta$ forces that $f(\beta) \in \dot{Q}_\beta$.  By elementarity of $j$, $N$ believes that $j(f) \restriction j(\beta) = f \restriction \beta$ is an element of the $j(\beta) = \beta$-th member of $j(\vec{\mathbb{P}})$, and that $j\big( f(\beta) \big) = f(\beta)$ is forced by $j(V_\beta \cap \mathbb{P}_\beta) = V_\beta \cap j(\vec{\mathbb{P}})_{j(\beta)} = V_\beta \cap j(\vec{\mathbb{P}})_{\beta} $ to be a member of $j(\dot{Q}_\beta)$.  This proves that $j(f) = f$ is a member of $j(\vec{\mathbb{P}})_\alpha$, so in particular a member of $j(\vec{\mathbb{P}})_\kappa$.

Finally, if we know that $\mathbb{P}_\kappa$ is $\kappa$-cc, then because $\text{crit}(j) = \kappa$ it follows that $j \restriction \mathbb{P}_\kappa = \text{id}_{\mathbb{P}_\kappa}$ is a regular embedding from $\mathbb{P}_\kappa \to j(\mathbb{P}_\kappa) = j(\vec{\mathbb{P}})_{j(\kappa)}$; i.e.\
\[
\mathbb{R}_0 :=j[\mathbb{P}_\kappa] \text{ is a regular suborder of }  \mathbb{R}_2:= j(\vec{\mathbb{P}})_{j(\kappa)}   .
\]
But because we now know that $\mathbb{R}_0 = j[\mathbb{P}_\kappa]$ is actually a suborder of $\mathbb{R}_1:= j(\vec{\mathbb{P}})_\kappa$ (by the previous part of the proof), and because $\mathbb{R}_1 = j(\vec{\mathbb{P}})_\kappa$ is regular in $\mathbb{R}_2 = j(\vec{\mathbb{P}})_{j(\kappa)}$ (because $j(\vec{\mathbb{P}})$ is a generalized iteration), it abstractly follows that $\mathbb{R}_0 = j[\mathbb{P}_\kappa]$ is a regular suborder of $\mathbb{R}_1 = j(\vec{\mathbb{P}})_\kappa$.\footnote{In general, if $\mathbb{R}_0$ is a suborder of $\mathbb{R}_1$, $\mathbb{R}_1$ is a suborder of $\mathbb{R}_2$, and $\mathbb{R}_0$ is a \emph{regular} suborder of $\mathbb{R}_2$, it follows easily that $\mathbb{R}_0$ is a regular suborder of $\mathbb{R}_1$.}   Since we showed earlier that $j[\mathbb{P}_\kappa] = V_\kappa \cap j(\vec{\mathbb{P}})_\kappa$, it follows that $\kappa$ is an active stage of $j(\vec{\mathbb{P}})$.
\end{proof}

The following lemma ensures preservation of closure properties:
\begin{lemma}\label{lem_LessMuClosedIteration}
Let $\mu < \kappa$ be regular cardinals.  Suppose $\langle \mathbb{P}_\alpha,\dot{Q}_\beta \ | \ \alpha \le \kappa, \beta < \kappa \rangle$ is a universal Kunen iteration with the following properties:
\begin{itemize}
 \item For all (active) $\alpha < \kappa$:
 \[
 \Vdash_{V_\alpha \cap \mathbb{P}_\alpha} \ \dot{Q}_\alpha \text{ is } < \mu \text{ closed}.
 \]
 \item The iteration uses $<\mu$ supports; i.e.\ for all limit $\alpha \le \kappa$, $\mathbb{P}_\alpha$ consists of those partial functions $f$ such that $\text{dom}(f)$ is a $<\mu$-sized subset of $\alpha$, and $f \restriction \bar{\alpha} \in \mathbb{P}_{\bar{\alpha}}$ for all $\bar{\alpha} \in \text{dom}(f)$.
\end{itemize}
Then $\mathbb{P}_\kappa$ is $<\mu$ closed.
\end{lemma}
\begin{proof}
This is basically the standard proof that $<\mu$-support iterations of $<\mu$-closed posets is $<\mu$-closed (as in section 2 of \cite{MR823775}), with a variation at successor steps due to the requirement of Definition \ref{def_KunenIteration} that we only use $V_\alpha \cap \mathbb{P}_\alpha$-names, rather than $\mathbb{P}_\alpha$-names. 

Fix $\lambda < \mu$ and suppose $\langle f_i \ | \ i < \lambda \rangle$ is a descending sequence in $\mathbb{P}_\kappa$.  Note that $\langle \text{dom}(f_i) \ | \ i < \lambda \rangle$ is a $\subseteq$-increasing sequence of $<\mu$-sized sets, and by regularity of $\mu$ it follows that
\[
d:= \bigcup_{i < \lambda} \text{dom}(f_i)
\]
has size $<\mu$.  We recursively define a function $f$ with domain $d$ and verify that $f$ is a condition in $\mathbb{P}_\kappa$ and is a lower bound for the $f_i$'s.  Suppose $\alpha \le \kappa$ and $f (  \bar{\alpha})$ is defined for all $\bar{\alpha} < \alpha$, and that the following induction hypothesis $\text{IH}_{\bar{\alpha}}$ holds for each $\bar{\alpha} < \alpha$:
\begin{align*}
\text{IH}_{\bar{\alpha}}: \ \  &  f \restriction (\bar{\alpha}+1) \text{ is a lower bound in } \mathbb{P}_{\bar{\alpha}+1} \text{ for the sequence }\langle f_i \restriction (\bar{\alpha}+1) \ | \ i < \lambda \rangle \\
& \text{ and } \text{dom}\big(f \restriction (\bar{\alpha}+1) \big) \subseteq d
\end{align*} 

First observe that the fact that $\text{IH}_{\bar{\alpha}}$ holds for all $\bar{\alpha} < \alpha$ immediately implies:
\begin{equation}\label{eq_f_cut_alpha}
f \restriction \alpha \text{ is a lower bound in } \mathbb{P}_\alpha \text{ for } \langle f_i \restriction \alpha \ | \ i < \lambda \rangle.
\end{equation}

If $\alpha = \kappa$ we are done, otherwise we need to define $f(\alpha)$ and verify that $\text{IH}_\alpha$ holds.  If $\alpha \notin d$ then $f(\alpha)$ isn't defined, i.e.\ $\alpha \notin \text{dom}(f)$.  Then $\text{IH}_\alpha$ is immediate.  So from now on assume $\alpha \in d$; then $\alpha \in \text{dom}(f_i)$ for all sufficiently large $i < \lambda$, say for all $i \ge i_{\alpha}$.  Note that by the definition of the ordering for $\mathbb{P}_{\alpha+1}$, together with the assumption that the $f_i$'s are descending, it follows that:
\begin{equation}\label{eq_DesendingSeq}
i_{\alpha} \le i < j < \lambda \ \ \implies f_j \restriction \alpha \Vdash^{\text{Convention}}_{V_{\alpha} \cap \mathbb{P}_{\alpha}}  f_j(\alpha) \le_{\dot{Q}_{\alpha}} f_i(\alpha).
\end{equation}
Then \eqref{eq_f_cut_alpha} and \eqref{eq_DesendingSeq} together imply
\[
f \restriction \alpha \Vdash^{\text{Convention}}_{V_{\alpha} \cap \mathbb{P}_{\alpha}} \ \langle f_i(\alpha) \ | \ i_{\alpha} \le i < \lambda \rangle \text{ is descending in } \dot{Q}_{\alpha}.
\]
Because $V_\alpha \cap \mathbb{P}_\alpha$ forces (by assumption) that $\dot{Q}_\alpha$ is $<\mu$ closed, then by the Maximality Principle there is a $V_\alpha \cap \mathbb{P}_\alpha$-name $\dot{q}$ such that
\[
f \restriction \alpha \Vdash^{\text{Convention}}_{V_{\alpha} \cap \mathbb{P}_{\alpha}} \  \dot{q} \text{ is a lower bound for }  \langle f_i(\alpha) \ | \ i_{\alpha} \le i < \lambda \rangle .
\]
Then set $f(\alpha):= \dot{q}$.  Then $\text{IH}_\alpha$ clearly holds.
\end{proof}

Here is how such universal iterations are typically used; we concentrate on Kunen's original construction for concreteness.  Suppose $j:V\to N$ is a huge embedding with critical point $\kappa$, and let $\delta = j(\kappa)$.  Kunen defined a finite support iteration $\langle \mathbb{P}_\alpha, \dot{Q}_\beta \ | \ \alpha \le \kappa, \beta < \kappa \rangle$ as in Definition \ref{def_KunenIteration}, where one starts by forcing with $\mathbb{P}_1=\text{Col}(\omega, < \kappa)$ and then at each active $\alpha$ one forces with $\mathbb{S}^{V^{V_\alpha \cap \mathbb{P}_\alpha}}(\alpha, < \kappa)$.  Kunen's proof that $\mathbb{P}_\kappa$ is $\kappa$-cc goes by induction; limit steps are taken care of by Theorem \ref{thm_FiniteSupport}, and successor steps $\alpha \to \alpha+1$ are taken care of by the induction hypothesis that $\mathbb{P}_\alpha$ is $\kappa$-cc, together with the fact that $\mathbb{S}^{V^{V_\alpha \cap \mathbb{P}_\alpha}}(\alpha, < \kappa)$ is $\kappa$-Knaster; so the product of $\mathbb{S}^{V^{V_\alpha \cap \mathbb{P}_\alpha}}(\alpha, < \kappa)$ with the $\kappa$-cc quotient $\mathbb{P}_\alpha/(V_\alpha \cap \mathbb{P}_\alpha)$ is $\kappa$-cc by Fact \ref{fact_KnasterProduct}.  Then by Lemma \ref{lem_SubsetVkappa}, $\kappa$ is an active stage of $j(\vec{\mathbb{P}})$, and so the identity regular embedding from $\mathbb{P}_\kappa \to j(\mathbb{P}_\kappa)$ lifts to a regular embedding from $\mathbb{P}_\kappa * \mathbb{S}^{V^{\mathbb{P}_\kappa}}(\kappa, < \delta) \to j(\mathbb{P}_\kappa)$.  Thus if $G'$ is generic over $N$ for $j(\mathbb{P}_\kappa)$ then in $N[G']$ there is a $G*H$ which is generic over $V$ for $\mathbb{P}_\kappa * \mathbb{S}^{V^{\mathbb{P}_\kappa}}(\kappa, < \delta)$; then $j$ lifts to $j_{G'}: V[G] \to N[G']$ and $N[G']$ sees that $m:= \cup j_{G'}[H]$ is a condition in its Silver collapse $\mathbb{S}^{N[G']}\big( \delta, < j(\delta) \big)$.  Forcing below $m$ yields a lifting of $j$ to a map from $V[G*H] \to N[G' * H']$, but one can approximate such a lifting in the model $V[G']$ (without actually forcing with $H'$).  Because this occurs in the $\delta$-cc extension $V[G']$, the resulting derived ideal in $V[G*H]$ will be $\delta = \kappa^{+V[G*H]}$-cc, i.e.\ saturated.  We refer the reader to Foreman~\cite{MattHandbook} for more details.

Here are some examples of such iterations in the literature on saturated ideals; the columns of the table correspond to the scheme described in Definition \ref{def_KunenIteration} above.  Kunen's original paper~\cite{MR495118} dealt only with ideals on $\omega_1$ and used finite support iteration; Laver showed how to do $<\mu$-support iterations for arbitrary regular $\mu$ in certain situations.  In all cases $\mu < \kappa$ are regular cardinals and $\mathbb{P}_0$ is the Levy collapse $\text{Col}(\mu, < \kappa)$ or the Silver collapse $\mathbb{S}(\mu, < \kappa)$ (see Section \ref{sec_Examples} for the definitions of these posets).  In the table, $\dot{g}_\alpha$ denotes the canonical $V_\alpha \cap \mathbb{P}_\alpha$-name for its generic object.

\vspace{15pt}

\begin{tabular}{p{.20\linewidth}|c|c}
Paper &  $\dot{\mathbb{Q}}_\alpha$ (when $V_\alpha \cap \mathbb{P}_\alpha \lessdot \mathbb{P}_\alpha$) & supports \\
\hline
\hline
Kunen~\cite{MR495118}  & $\mathbb{S}^{V[\dot{g}_\alpha]}(\alpha, < \kappa)$ & finite (i.e.\ $<\mu = \omega$) \\
\hline
Magidor (see \cite{MattHandbook})  & $\text{Col}^{V[\dot{g}_\alpha]}(\alpha, < \kappa)$ & $<\mu$ \\
\hline
Laver~\cite{MR673792}   & $\mathbb{E}^{V[\dot{g}_\alpha]}(\alpha, < \kappa)$ & $<\mu$  \\
\hline
Foreman-Komjath~\cite{MR2151585}  & $\mathbb{S}^{V[\dot{g}_\alpha]}(\alpha, < \kappa) \times \mathbb{B}^{V[\dot{g}_\alpha]}(\mu, \alpha^{<\alpha}, \kappa)$ & $<\mu$ \\
\hline
Foreman~\cite{MR730584} and Foreman-Laver~\cite{MR925267}   &  A certain product of Silver collapses  & $<\mu$ \\
\hline
\end{tabular}

\vspace{15pt}

The following is the key technical theorem of this section:
\begin{theorem}\label{thm_MainTheorem}
Let $\kappa$ be Mahlo and $\langle \mathbb{P}_\alpha, \dot{\mathbb{Q}}_\beta \ | \ \alpha \le \kappa, \ \beta < \kappa \rangle$ be a universal Kunen iteration as in Definition \ref{def_KunenIteration}, which uses some mix of inverse and direct limits.  Assume $W \prec (H_{\kappa^+}, \in, \vec{\mathbb{P}})$ is such that $|W|=W \cap \kappa =: \gamma_W \in \kappa$ and ${}^{<|W|} W \subset W$.  Assume that:
\begin{enumerate}
 \item $\mathbb{P}_\kappa$ is a direct limit;
 \item $\mathbb{P}_{\gamma_W}$ is a direct limit;
 \item\label{item_ConditionsInVkappa} for every $\alpha < \gamma_W$:
\begin{equation*}
\Vdash_{V_\alpha \cap \mathbb{P}_\alpha} \dot{Q}_\alpha \subseteq V_\kappa[\dot{g}_\alpha]
\end{equation*}
where $\dot{g}_\alpha$ is the canonical $V_\alpha \cap \mathbb{P}_\alpha$-name for its generic;

  \item\label{item_ReductionInNextStep} for every $\alpha < \gamma_W$:
\begin{equation*}
\Vdash_{V_\alpha\cap \mathbb{P}_\alpha} \check{W}[\dot{g}_\alpha] \cap \dot{Q}_\alpha \text{ is a regular suborder of } \dot{Q}_\alpha.   
\end{equation*}
(Note this is trivially true for non-active $\alpha$).
\end{enumerate}

Then $W \cap \mathbb{P}_\kappa$ is a regular suborder of $\mathbb{P}_\kappa$.  
\end{theorem}

\begin{proof}
Note that Lemma \ref{lem_SubsetVkappa} implies $\mathbb{P}_\kappa \subset V_\kappa$.  Fix any $W \prec (H_{\kappa^+}, \in, \mathbb{P}_\kappa)$ as in the hypotheses of the theorem, and set $\gamma_W:= W \cap \kappa$.  We will  recursively construct a coherent, conservative system of reduction operations
\begin{equation}\label{eq_SystemToConstruct}
\langle \pi^W_\alpha: \mathbb{P}_\alpha \to W \cap \mathbb{P}_\alpha \ | \alpha < \gamma_W \rangle.
\end{equation}
At the end of the proof we will use the system from \eqref{eq_SystemToConstruct} to get a reduction operation from $\mathbb{P}_\kappa \to W \cap \mathbb{P}_\kappa$.

Let $\beta < \gamma_W$ and assume that the system \eqref{eq_SystemToConstruct} has been constructed below $\beta$; we need to define $\pi^W_\beta$.  Assume first that $\beta$ is a successor ordinal, say $\beta = \alpha+1$.  So $\pi^W_\alpha: \mathbb{P}_\alpha \to W \cap \mathbb{P}_\alpha$ is a reduction operation.  If $\alpha$ is passive then $\mathbb{P}_\beta = \mathbb{P}_\alpha$ and $\pi^W_\beta:= \pi^W_\alpha$ works.  Otherwise $\alpha$ is active and assumption \ref{item_ReductionInNextStep} and the Maximality Principle imply there is some $V_\alpha \cap \mathbb{P}_\alpha$-name $\dot{\rho}$ such that
\begin{equation*}
\Vdash_{V_\alpha \cap \mathbb{P}_\alpha} \dot{\rho}: \dot{Q}_\alpha \to \check{W}[\dot{g}_\alpha] \cap \dot{Q}_\alpha \text{ is a reduction operation}.
\end{equation*}
Given $f \in \mathbb{P}_{\alpha+1}$, define $\pi^W_{\alpha+1}(f)$ as follows.  The first $\alpha$ coordinates are given by $\pi^W_\alpha(f \restriction \alpha)$.  Now consider the $V_\alpha \cap \mathbb{P}_\alpha$-name $\dot{\rho}\big( f(\alpha) \big)$; assumption \ref{item_ConditionsInVkappa} implies $\dot{\rho}\big( f(\alpha) \big)$ is an element of $V_\kappa[\dot{g}_\alpha] \cap \check{W}[\dot{g}_\alpha]$.  Because $\alpha \in W \in \Gamma$ then $V_\alpha \cap \mathbb{P}_\alpha \subset V_{\gamma_W} = W \cap V_\kappa$; so $V_\kappa[\dot{g}_\alpha] \cap \check{W}[\dot{g}_\alpha]$ is forced to be the same as $(W \cap V_\kappa)[\dot{g}_\alpha] = V_{\gamma_W}[\dot{g}_\alpha]$.  By Fact \ref{fact_NamesInModel}, there is some $V_\alpha \cap \mathbb{P}_\alpha$-name $\dot{y}^W_{f(\alpha)} \in V_{\gamma_W} \subset W$ such that 
\begin{equation}\label{eq_NamesAreEqual}
\Vdash_{V_\alpha \cap \mathbb{P}_\alpha} \dot{y}^W_{f(\alpha)} = \dot{\rho}\big( f(\alpha) \big).
\end{equation}
We use this $\dot{y}^W_{f(\alpha)}$ as the $\alpha$-th coordinate; i.e.\
\begin{equation*}
\pi^W_{\alpha+1}(f) := \pi^W_{\alpha}\big( f \restriction \alpha \big) \cup \{ (\alpha, \dot{y}^W_{f(\alpha)} ) \}.
\end{equation*}

Note that because $\dot{y}^W_{f(\alpha)} \in W$ and $\pi^W_\alpha$ maps into $W \cap \mathbb{P}_\alpha$ by the induction assumption, then $\pi^W_{\alpha+1}$ maps into $W \cap \mathbb{P}_{\alpha+1}$.  We must verify that $\vec{\pi}^W \restriction (\alpha+1) \cup \{ (\alpha+1, \pi^W_{\alpha+1} )\} $ is conservative and coherent.

Coherence of $\vec{\pi}^W \restriction (\alpha+1) \cup \{ (\alpha+1, \pi^W_{\alpha+1}) \} $ follows immediately from the inductively-assumed coherency of $\vec{\pi}^W \restriction (\alpha+1)$ together with the definition of $\pi^W_{\alpha+1}$.  Finally we check conservativity of $\vec{\pi}^W \restriction (\alpha+1) \cup \{ (\alpha+1, \pi^W_{\alpha+1} )\} $ (which will also show that $\pi^W_{\alpha+1}$ is a reduction operation, as defined immediately before Definition \ref{def_ConservCoherent}).  Assume:
\begin{itemize}
 \item $f \in \mathbb{P}_{\alpha+1}$;
 \item $h \le_{\mathbb{P}_{\alpha+1}} \pi^W_{\alpha+1}(f)$, where $h \in W$;
 \item $\eta \le \alpha$, $f'_\eta \le_{\mathbb{P}_\eta} f \restriction \eta$, and $f'_\eta \le_{\mathbb{P}_\eta} h \restriction \eta$
\end{itemize}
We need to find some $f' \in \mathbb{P}_{\alpha+1}$ such that $f' \le_{\mathbb{P}_{\alpha+1}} f$, $f' \restriction \eta = f'_\eta$, and $f' \le_{\mathbb{P}_{\alpha+1}} h$.  

Our inductive hypothesis that $\vec{\pi}^W \restriction (\alpha+1)$ is conservative ensures that there is an $f'_\alpha \in \mathbb{P}_\alpha$ such that $f'_\alpha \le_{\mathbb{P}_\alpha} f \restriction \alpha$, $f'_\alpha \restriction \eta = f'_\eta$, and $f'_\alpha \le_{\mathbb{P}_\alpha} h \restriction \alpha$. 

Now consider the $V_\alpha \cap \mathbb{P}_\alpha$-names $f(\alpha)$ and $h(\alpha)$.  Now $h \restriction \alpha \le_{\mathbb{P}_\alpha} \pi^W_\alpha(f \restriction \alpha)$ and 
\begin{equation}\label{eq_h_restrict_alpha_forces}
h \restriction \alpha \Vdash_{\mathbb{P}_\alpha} h(\alpha) \le_{\dot{Q}_\alpha} \pi^W_{\alpha+1} (f)  (\alpha) = \dot{y}^W_{f(\alpha)}
\end{equation}
where again we are using Convention \ref{conv_ViewNames}.  Let $\text{WeakestReduct}_{\text{ro}(V_\alpha \cap \mathbb{P}_\alpha)} \big( h \restriction \alpha \big)$ be the maximum reduct of $h \restriction \alpha$ into $\text{ro}(V_\alpha \cap \mathbb{P}_\alpha)$ as in Definition \ref{def_WeakestReduct}.  Then \eqref{eq_h_restrict_alpha_forces} and Lemma \ref{lem_IfWeakestReductForcesThenOrigDoes} imply that
\begin{equation}\label{eq_BelowDotRhofalpha}
\text{WeakestReduct}_{\text{ro}(V_\alpha \cap \mathbb{P}_\alpha)} \big( h \restriction \alpha \big) \Vdash_{\text{ro}(V_\alpha \cap \mathbb{P}_\alpha)} h(\alpha) \le_{\dot{Q}_\alpha} \pi^W_{\alpha+1} (f)  (\alpha) = \dot{y}^W_{f(\alpha)} = \dot{\rho}\big( f(\alpha) \big).
\end{equation}
Because $h \in W$ and $\alpha \in W$, it follows that $h(\alpha) \in W$, which in turn implies 
\begin{equation}\label{eq_InW}
\text{WeakestReduct}_{\text{ro}(V_\alpha \cap \mathbb{P}_\alpha)} \Vdash_{\text{ro}(V_\alpha \cap \mathbb{P}_\alpha)} h(\alpha) \in \dot{Q}_\alpha \cap \check{W}[\dot{g}_\alpha].
\end{equation}
Because $\dot{\rho}$ is forced by $V_\alpha \cap \mathbb{P}_\alpha$ to be a reduction from $\dot{Q}_\alpha \to \dot{Q}_\alpha \cap W[\dot{g}_\alpha]$, then \eqref{eq_BelowDotRhofalpha} and \eqref{eq_InW} imply
\begin{equation}
\text{WeakestReduct}_{\text{ro}(V_\alpha \cap \mathbb{P}_\alpha)} \big( h \restriction \alpha \big) \Vdash_{\text{ro}(V_\alpha \cap \mathbb{P}_\alpha)} \exists q \le_{\dot{Q}_\alpha} f(\alpha) \text{ such that } q \le_{\dot{Q}_\alpha} h(\alpha).
\end{equation}
By the Maximality Principle there is a $V_\alpha \cap \mathbb{P}_\alpha$-name $\dot{q}_{h(\alpha), f(\alpha)}$ such that
\begin{equation}\label{eq_WeakestReductForces}
\text{WeakestReduct}_{\text{ro}(V_\alpha \cap \mathbb{P}_\alpha)} \big( h \restriction \alpha \big) \Vdash_{\text{ro}(V_\alpha \cap \mathbb{P}_\alpha)} \dot{q}_{h(\alpha), f(\alpha)} \le_{\dot{Q}_\alpha} f(\alpha) \text{ and } \dot{q}_{h(\alpha), f(\alpha)} \le_{\dot{Q}_\alpha} h(\alpha)
\end{equation}
and by standard `definition of names by cases" we can WLOG assume that $\emptyset \Vdash_{V_\alpha \cap \mathbb{P}_\alpha} \dot{q}_{h(\alpha), f(\alpha)} \in \dot{Q}_\alpha$ and that $\emptyset \Vdash_{V_\alpha \cap \mathbb{P}_\alpha} \dot{q}_{h(\alpha), f(\alpha)} \le_{\dot{Q}_\alpha} f(\alpha)$.  So
\begin{equation*}
f':= f'_\alpha \cup \{ \big(\alpha, \dot{q}_{h(\alpha), f(\alpha)} \big)  \}
\end{equation*}
is a condition in $\mathbb{P}_{\alpha+1}$, $f' \le_{\mathbb{P}_{\alpha+1}} f$, and $f' \restriction \eta = f'_\alpha \restriction \eta =   f'_\eta$.  We only have left to check that
\begin{equation}\label{eq_InequalityToShow}
f' \le_{\mathbb{P}_{\alpha+1}} h.
\end{equation}
Now 
\begin{equation}\label{eq_Below_h_restrict_alpha}
f' \restriction \alpha = f'_\alpha \le_{\mathbb{P}_\alpha} h \restriction \alpha
\end{equation}
and $f'(\alpha) = \dot{q}_{h(\alpha), f(\alpha)}$; so we only have left to prove that
\begin{equation}
f'_\alpha \Vdash_{\mathbb{P}_\alpha}^{\text{Convention}} \dot{q}_{h(\alpha), f(\alpha)}   \le_{\dot{Q}_\alpha} h(\alpha).
\end{equation}
Moreover, by \eqref{eq_Below_h_restrict_alpha} it in turn suffices to prove
\begin{equation}
h \restriction \alpha \Vdash_{\mathbb{P}_\alpha}^{\text{Convention}}  \dot{q}_{h(\alpha), f(\alpha)}  \le_{\dot{Q}_\alpha} h(\alpha).
\end{equation}
But this follows from \eqref{eq_WeakestReductForces} and Lemma \ref{lem_IfWeakestReductForcesThenOrigDoes}.  This completes the proof that $\vec{\pi}^W \restriction (\alpha+1) \cup \{  (\alpha+1, \pi^W_{\alpha+1}) \}$ is a coherent conservative system of reduction operations.  This completes the case where $\beta$ is a successor ordinal.

Now assume $\beta$ is a limit ordinal (still assuming $\beta < \gamma_W$).  Because stage $\beta$ of $\mathbb{P}_\beta$ is either an inverse or a direct limit, and because ${}^{<\gamma_W} W \subset W$ and $\beta \in W$, it follows that $\mathbb{P}_\beta$ and $\mathbb{P}_\beta \cap W$ use the same limit scheme (with respect to $\langle \mathbb{P}_\alpha \ | \ \alpha < \beta \rangle$ and $\langle \mathbb{P}_\alpha \cap W \ | \ \alpha < \beta \rangle$, respectively).  So Lemma \ref{lem_WhatToDoAtLimits} applies, and ensures that the natural limit of $\langle \pi^W_\alpha \ | \ \alpha < \beta \rangle$ (depending on the limit scheme used at $\beta$) is a reduction, and thus we obtain a coherent, conservative system of reductions $\langle \pi^W_\alpha \ | \ \alpha \le \beta \rangle$.

We have constructed a coherent conservative system $\langle \pi^W_\alpha \ | \ \alpha < \gamma_W \rangle$.  Now by the assumptions that ${}^{<\gamma_W} W \subset W$ and that $\mathbb{P}_{\gamma_W}$ is a direct limit, then $\mathbb{P}_{\gamma_W}$ and $W \cap \mathbb{P}_{\gamma_W}$ are direct limits of $\langle \mathbb{P}_\alpha \ | \ \alpha < \gamma_W \rangle$ and $\langle W \cap \mathbb{P}_\alpha \ | \ \alpha < \gamma_W \rangle$, respectively.  So we may again invoke Lemma \ref{lem_WhatToDoAtLimits} to extend this to a coherent conservative system $\langle \pi^W_\alpha \ | \ \alpha \le \gamma_W \rangle$.  Although $\mathbb{P}_{\gamma_W} \notin W$, because $\mathbb{P}_{\gamma_W}$ is a direct limit then we still have that $\text{id}: W \cap \mathbb{P}_{\gamma_W} \to \mathbb{P}_{\gamma_W}$ is order and incompatibility preserving.  This, together with the existence of the reduction operation $\pi^W_{\gamma_W}: \mathbb{P}_{\gamma_W} \to W \cap \mathbb{P}_{\gamma_W}$, implies that 
\begin{equation*}
W \cap \mathbb{P}_{\gamma_W} \text{ is a regular suborder of } \mathbb{P}_{\gamma_W}.
\end{equation*}
Because $\mathbb{P}_{\gamma_W}$ is a regular suborder of $\mathbb{P}_\kappa$ it follows that
\begin{equation}\label{eq_WCapGammaWRegInPkappa}
W \cap \mathbb{P}_{\gamma_W} \text{ is a regular suborder of } \mathbb{P}_\kappa.
\end{equation}

Finally, observe that because both $\mathbb{P}_{\gamma_W}=\mathbb{P}_{W \cap \kappa}$ and $\mathbb{P}_\kappa$ are direct limits, then
\begin{equation}\label{eq_WCapPkappaEqualGammaW}
W \cap \mathbb{P}_\kappa = W \cap \mathbb{P}_{\gamma_W}.
\end{equation}
Note this is a literal equality, because our conditions are partial functions (i.e.\ direct limits are simply unions).\footnote{If $f \in W \cap \mathbb{P}_\kappa$ then $W$ believes there is some $\beta < \kappa$ such that $f \in \mathbb{P}_\beta$; so this $\beta$ must be strictly below $W \cap \kappa = \gamma_W$.  Thus $f \in W \cap \mathbb{P}_{\gamma_W}$.}  So \eqref{eq_WCapGammaWRegInPkappa} and \eqref{eq_WCapPkappaEqualGammaW} together imply that $W \cap \mathbb{P}_\kappa$ is a regular suborder of $\mathbb{P}_\kappa$.
\end{proof}

The following theorem is almost an immediate corollary of Theorem \ref{thm_MainTheorem}; however it is typically less useful than Theorem \ref{thm_WeakCompact_Kunen_MainThm_Precise}, because typically $\kappa$ is a very large cardinal in the context of universal Kunen iterations.

\begin{theorem}\label{thm_MainTheorem_Mahlo_Precise}
Assume $\kappa$ is Mahlo, $\langle \mathbb{P}_\alpha, \dot{Q}_\beta \ | \ \alpha \le \kappa, \beta < \kappa \rangle$ is a universal Kunen iteration as in Definition \ref{def_KunenIteration} which uses some mix of inverse and direct limits, and that $\mathbb{P}_\gamma$ is a direct limit for all inaccessible $\gamma \le \kappa$.  

Let $\Phi(W,\kappa)$ denote the statement that $W \prec H_{\kappa^+}$, $|W|=W \cap \kappa \in \kappa$, and ${}^{<|W|} W \subset W$.  Assume that for each $\alpha < \kappa$:
\begin{equation*}
\Vdash_{V_\alpha \cap \mathbb{P}_\alpha} \dot{Q}_\alpha \text{ is layered a.e. on } \{ M  \ | \  \Phi^{V[\dot{g}_\alpha]}(M,\kappa) \} 
\end{equation*}
and
\begin{equation*}
\Vdash_{V_\alpha \cap \mathbb{P}_\alpha} \dot{Q}_\alpha \subset V_\kappa[\dot{g}_\alpha]
\end{equation*}
where $\dot{g}_\alpha$ is the canonical $V_\alpha \cap \mathbb{P}_\alpha$-name for its generic.

Then $\mathbb{P}_\kappa$ is layered a.e.\ on the stationary set
\begin{equation*}
\Gamma:= \{  W \ | \ \Phi^V(W,\kappa) \}.
\end{equation*}
i.e.\ for all but nonstationarily many $W \in \Gamma$, $W \cap \mathbb{P}_\kappa$ is a regular suborder of $\mathbb{P}_\kappa$.  In particular, $\mathbb{P}_\kappa$ is $\kappa$-Knaster by Lemma \ref{lem_LayeredImpliesKnaster}.
\end{theorem}
\begin{proof}
Note that $\Gamma$ is stationary because $\kappa$ is Mahlo, and the assumptions guarantee (by Lemma \ref{lem_SubsetVkappa}) that $\mathbb{P}_\kappa \subset V_\kappa$.  For each active $\alpha < \kappa$ let $\dot{\mathfrak{A}}_\alpha$ be a $V_\alpha \cap \mathbb{P}_\alpha$-name for a first order structure on $H_{\kappa^+}[\dot{g}_\alpha]$ witnessing the ``a.e." part of the assumption about $\dot{Q}_\alpha$; i.e.\ so that
\begin{equation}\label{eq_ae}
\Vdash_{V_\alpha \cap \mathbb{P}_\alpha} \forall M  \ \ M \prec \dot{\mathfrak{A}}_\alpha \wedge \Phi(M,\kappa) \implies \dot{Q}_\alpha \cap M \text{ is a regular suborder of } \dot{Q}_\alpha.
\end{equation}

Now let $W \in \Gamma$ such that $W \prec (H_{\kappa^+}, \in, \vec{\mathbb{P}}, \vec{\dot{\mathfrak{A}}})$.  Let $\alpha \in W \cap \kappa =: \gamma_W$.  Because $|V_\alpha \cap \mathbb{P}_\alpha| < \gamma_W$ then in particular $V_\alpha \cap \mathbb{P}_\alpha$ has the $\gamma_W$-cc, so together with $<\gamma_W$ closure of $W$ and the fact that $\dot{\mathfrak{A}}_\alpha \in W$ we have:
\begin{equation*}
\Vdash_{V_\alpha \cap \mathbb{P}_\alpha}  \Phi(W[\dot{g}_\alpha], \kappa) \text{ and } W[\dot{g}_\alpha] \prec \dot{\mathfrak{A}}_\alpha.
\end{equation*}
Thus by \eqref{eq_ae} we obtain:
\begin{equation}\label{eq_alphaLessRegular}
\forall \alpha < \gamma_W \ \ \  \Vdash_{V_\alpha \cap \mathbb{P}_\alpha} \dot{Q}_\alpha \cap W[\dot{g}_\alpha] \text{ is a regular suborder of } \dot{Q}_\alpha.
\end{equation}
Then \eqref{eq_alphaLessRegular}, together with our assumptions that $\mathbb{P}_\gamma$ is a direct limit for all inaccessible $\gamma \le \kappa$, implies by Theorem \ref{thm_MainTheorem} that $W \cap \mathbb{P}_\kappa$ is a regular suborder of $\mathbb{P}_\kappa$.
\end{proof}

We finally state and prove the precise version of Theorem \ref{thm_WeakCompact_Kunen_MainThm_Imprecise}:
It tells us that if $\kappa$ is weakly compact and direct limits are used at all inaccessibles, then \textbf{any} universal Kunen iteration of $\kappa$-cc posets will be $\kappa$-stationarily layered (and thus $\kappa$-Knaster).  We emphasize that assumption \ref{item_KappaCCForcedByInner} in the statement of Theorem \ref{thm_WeakCompact_Kunen_MainThm_Precise} is only assuming that $V_\alpha \cap \mathbb{P}_\alpha$ forces $\dot{Q}_\alpha$ to be $\kappa$-cc; it is NOT assuming that $\dot{Q}_\alpha$ remains $\kappa$-cc in $V^{\mathbb{P}_\alpha}$ (though this will be the case automatically in the end).

\begin{theorem}\label{thm_WeakCompact_Kunen_MainThm_Precise}
Suppose $\kappa$ is weakly compact and $\langle \mathbb{P}_\alpha, \dot{Q}_\beta \ | \ \alpha \le \kappa, \beta < \kappa \rangle$ is a universal Kunen iteration as in Definition \ref{def_KunenIteration} which uses some mix of inverse and direct limits.  Suppose
\begin{enumerate}
 \item direct limits are taken at all inaccessible $\gamma \le \kappa$;
 \item for every active $\alpha < \kappa$:  $\Vdash_{V_\alpha \cap \mathbb{P}_\alpha} \dot{Q}_\alpha \subset V_\kappa[\dot{g}_\alpha]$;
 \item\label{item_KappaCCForcedByInner} Each $\dot{Q}_\alpha$ is forced by $V_\alpha \cap \mathbb{P}_\alpha$ to be $\kappa$-cc.\footnote{Notice again this is a weaker assumption than requiring $\mathbb{P}_\alpha$ to force that $\dot{Q}_\alpha$ is $\kappa$-cc, though in the end that will be true as well.}
\end{enumerate}
Then $\mathbb{P}_\kappa \subset V_\kappa$ and is layered on some stationary subset of
\begin{equation*}
\Gamma:= \{ W \in P_\kappa(V_\kappa) \ | \ W = V_\gamma \text{ for some inaccessible } \gamma < \kappa   \}.
\end{equation*}
In particular $\mathbb{P}_\kappa$ is $\kappa$-Knaster by Lemma \ref{lem_LayeredImpliesKnaster}.
\end{theorem}
\begin{proof}
First note that $\Gamma$ is stationary, because $\kappa$ is Mahlo.  Also Lemma \ref{lem_SubsetVkappa} implies that $\mathbb{P}_\kappa \subset V_\kappa$.  Suppose toward a contradiction that there were some algebra $\mathfrak{A} = (V_\kappa, \in, \dots)$ such that 
\begin{equation}
\forall W \in \Gamma \ W \prec \mathfrak{A} \implies W \cap \mathbb{P}_\kappa \text{ is NOT a regular suborder of } \mathbb{P}_\kappa.
\end{equation}    

Fix some transitive $\kappa$-sized, $<\kappa$-closed models $H$ and $H'$ such that:
\begin{itemize}
 \item $\kappa \subset H \in H'$;
 \item $H \prec H' \prec (H_{\kappa^+}, \in, \{ \mathfrak{A}, \vec{\mathbb{P}} \})$
\end{itemize}

Because $\kappa$ is weakly compact, there is some transitive $N'$ and an elementary $j: H' \to N'$ with critical point $\kappa$.  Because $H \in H'$, $\text{crit}(j) = \kappa$, and $H' \models |H|=\kappa$, then
\begin{equation}
j[H] \in N' \ .
\end{equation}
Note also that $P_\kappa(H) \in H'$, and that $P_\kappa(H) \subset H'$ because $H'$ is $<\kappa$-closed and $H \subset H'$.   Define
\begin{equation}
U:= \{ X \in H' \cap P\big(P_\kappa(H)\big) \ | \  j[H] \in j(X)   \}.
\end{equation}

Then $U$ is an ultrafilter on $H' \cap P\big(P_\kappa(H)\big)$ which is $<\kappa$-closed in $V$, and normal with respect to $\kappa$ sequences from $H'$; i.e.\ if $F \in H'$ is a regressive function on some $Z \in U$, then there is a $Z' \in U$ on which $F$ is constant.\footnote{It is not really necessary to work with both $H$ and $H'$ here.  We could really use $V_\kappa$ in place of the $H$, and work with the ultrafilter on $H' \cap P(\kappa)$ derived from $j$ (rather than the ultrafilter on $H' \cap P(P_\kappa(H))$ derived from $j$), then make an argument at the end involving Skolem hulls of $\gamma < \kappa$.  We use $H$ merely so that the $U$-measure one sets concentrate directly on the kinds of model described in Theorem \ref{thm_MainTheorem}, so we can directly apply that theorem after the following claim.}

\begin{nonGlobalClaim}\label{clm_U_many}
Let $Z$ be the set of $W \in H' \cap P_\kappa(H)$ which satisfy clause \ref{item_ReductionInNextStep} of the assumptions of Theorem \ref{thm_MainTheorem}.  Then $Z \in U$.
\end{nonGlobalClaim}
\begin{proof}
(of Claim \ref{clm_U_many})  First note that $Z$ is an element of $H'$, because $\vec{\mathbb{P}} \in H'$.  Suppose for a contradiction that $Z \notin U$.  Then because $Z \in H'$ and $U$ is an ultrafilter on $H' \cap P\big(P_\kappa(H)\big)$, there are $U$-many $W \in \Gamma$ such that clause \ref{item_ReductionInNextStep} fails for some $\alpha < W \cap \kappa$; that is, $U$-many $W$ such that
\begin{multline*}
\exists \alpha_W < \gamma_W \ \exists p_W \in V_{\alpha_W} \cap \mathbb{P}_{\alpha_W}  \\ p_W \Vdash_{V_{\alpha_W} \cap \mathbb{P}_{\alpha_W}} \dot{Q}_{\alpha_W} \cap W[\dot{g}_{\alpha_W}] \text{ is NOT a regular suborder of } \dot{Q}_{\alpha_W}.
\end{multline*}

Notice that the map $W \mapsto (\alpha_W, p_W)$ is regressive, and also an element of $H'$.  So by normality of $U$ with respect to regressive functions from $H'$, there is some fixed $\alpha^*$ and a fixed $p^* \in V_{\alpha^*} \cap \mathbb{P}_{\alpha^*}$ such that for $U$-many $W \in \Gamma$:
\begin{equation*}
p^* \Vdash_{V_{\alpha^*} \cap \mathbb{P}_{\alpha^*}} \dot{Q}_{\alpha^*} \cap W[\dot{g}_{\alpha^*}] \text{ is NOT a regular suborder of } \dot{Q}_{\alpha^*}.
\end{equation*}
Let $A$ denote the set of such $W$ (so $A \in U$, i.e.\ $j[H] \in j(A)$).  Notice that $A \in H'$ (it is definable from the regressive function and $(\alpha^*, p^*)$).

Let $g^*$ be generic over $V$ for $V_{\alpha^*} \cap \mathbb{P}_{\alpha^*}$ with $p^* \in g^*$.  Because $p^* \in g^*$ it follows that
\begin{equation}\label{eq_ToContradict}
\forall W \in A \ \ \mathbb{Q}_{\alpha^*} \cap W[g^*] \text{ is NOT a regular suborder of } \mathbb{Q}_{\alpha^*}
\end{equation}
where $\mathbb{Q}_{\alpha^*}$ is the evaluation of $\dot{Q}_{\alpha^*}$ by $g^*$.

Because $V_{\alpha^*} \cap \mathbb{P}_{\alpha^*} \in V_\kappa$ then $j$ lifts in $V[g^*]$ to an elementary
\begin{equation*}
\tilde{j}: H'[g^*] \to N'[g^*] = N'[\tilde{j}(g^*)].
\end{equation*}

Because $H$ is $<\kappa$-closed in $V$ and $g^*$ is generic for a $\kappa$-cc poset (in fact a poset of size $<\kappa$), then $H[g^*]$ is $<\kappa$ closed in $V[g^*]$, and in particular in $N'[g^*]$.  Since also $\mathbb{Q}_{\alpha^*}$ is $\kappa$-cc in $V[g^*]$, and thus in $N'[g^*]$, then
\begin{equation}
N'[g^*] \text{ and }  H[g^*] \text{ have the same maximal antichains of } \mathbb{Q}_{\alpha^*}. 
\end{equation}

Then $H[g^*]$, $\tilde{j} \restriction H[g^*]$, $N'[g^*]$, and $\mathbb{Q}_{\alpha^*}$ satisfy the assumptions of Lemma \ref{lem_EmbeddingsAndKappaCC}; note we are using here that $j[H] \in N'$ to conclude that $j[H][g^*] = \tilde{j}(H[g^*]) \in N'[g^*]$.  So by Lemma \ref{lem_EmbeddingsAndKappaCC}:
\begin{equation}\label{eq_TildejHgStar}
\tilde{j}[H][g^*] \cap \tilde{j}(\mathbb{Q}_{\alpha^*}) \text{ is a regular suborder of } \tilde{j}(\mathbb{Q}_{\alpha^*}).
\end{equation}

Because $\tilde{j}[H] = j[H] \in j(A) =  \tilde{j}(A)$ and $\tilde{j}(g^*)=g^*$, then by \eqref{eq_TildejHgStar} we get:
\begin{equation*}
N[g^*] \models  \ \exists X \in \tilde{j}(A) \ \ X[\tilde{j}(g^*)] \cap \tilde{j}(\mathbb{Q}_{\alpha^*}) \text{ is a regular suborder of } \tilde{j}(\mathbb{Q}_{\alpha^*}).
\end{equation*}
By elementarity of $\tilde{j}$ then
\begin{equation*}
H'[g^*] \models \ \exists X \in A \ \ X[g^*] \cap \mathbb{Q}_{\alpha^*} \text{ is a regular suborder of } \mathbb{Q}_{\alpha^*}
\end{equation*}
which contradicts \eqref{eq_ToContradict} and completes the proof of Claim \ref{clm_U_many}.
\end{proof}

Let $B$ be the $U$-set given by Claim \ref{clm_U_many}.  Because $B \in U$ then $H' \models$ ``$B$ is stationary in $P_\kappa(H)$" and because $H' \prec H_{\kappa^+}$ then $B$ really is stationary in $P_\kappa(H)$.  By Theorem \ref{thm_MainTheorem}, $\mathbb{P}_\kappa$ is layered a.e.\ on $B$.  Thus $\mathbb{P}$ is $\kappa$ stationarily layered.   
\end{proof}

\begin{remark}
If the $\kappa$ in the assumptions of Theorem \ref{thm_WeakCompact_Kunen_MainThm_Precise} is measurable, and $U$ is any normal ultrafilter on $\kappa$, then essentially the same proof shows that $\mathbb{P}_\kappa$ is layered on a $U$-measure one set; i.e.\ there are $U$-many $\alpha < \kappa$ such that $V_\alpha \cap \mathbb{P}_\kappa$ is a regular suborder of $\mathbb{P}_\kappa$ (i.e.\ $U$-many active $\alpha$).  
\end{remark}

\section{Applications}\label{sec_Apps}

In this section we use our main iteration result (Theorem \ref{thm_WeakCompact_Kunen_MainThm_Precise}) to prove Theorems \ref{thm_IdealAbsorb} and \ref{thm_ChangTree}.  A few remarks are in order first:  
\begin{enumerate}
 \item The special case where $\mu = \omega$ in the conclusion of Theorem \ref{thm_IdealAbsorb} follows easily from Kunen's original argument (\cite{MR495118}), without having to use our new Theorem \ref{thm_WeakCompact_Kunen_MainThm_Precise}.   The reason that the case $\mu = \omega$ is significantly easier is because of the Solovay-Tennenbaum Theorem \ref{thm_FiniteSupport} about finite support iterations.  
 \item Similarly, the special case $\mu = \omega$ of Theorem \ref{thm_ChangTree} could have easily been proven using Kunen's original construction, together with Kanamori's Theorem \ref{thm_Kanamori} below, without having to use our new Theorem \ref{thm_WeakCompact_Kunen_MainThm_Precise}.  The reason, again, is that in this particular case (i.e.\ when $\kappa$ becomes $\omega_1$), finite support iterations are used and so the Solovay-Tennenbaum Theorem \ref{thm_FiniteSupport} is applicable.
 \item Certain specific instances of Theorem \ref{thm_IdealAbsorb} were already known for arbitrary $\mu$; e.g.\ such as when the poset $\mathbb{S}^{H_{\kappa^{++}}}_{r,\phi}$ from the statement of Theorem \ref{thm_IdealAbsorb} is the poset $\text{Col}(\mu, < \kappa^+)$\footnote{Equivalently, the $<\mu$-support product of $\kappa^+$ many copies of $\text{Col}(\mu,\kappa)$} and certain other collapsing-type posets.  The examples listed in Corollary \ref{cor_Examples} are all new, however, and the template provided by Theorem \ref{thm_IdealAbsorb} for absorbing posets into quotients of saturated ideals is very general.
\end{enumerate}

\subsection{Proof of Theorem \ref{thm_IdealAbsorb}}

Fix $\phi(-,-)$ as in the statement of the theorem.  Assume $j: V \to N$ is a huge embedding with critical point $\kappa$.  Let $\delta = j(\kappa)$.  Assume $\mu$ is a regular cardinal and $\mu < \kappa$.  Define a $<\mu$-support universal Kunen iteration  
\[
\langle \mathbb{P}_\alpha, \dot{Q}_\beta \ | \ \alpha \le \kappa, \beta < \kappa \rangle
\]
(as in Definition \ref{def_KunenIteration}) as follows:  $\mathbb{P}_0$ is $\text{Col}(\mu, < \kappa)$.  For active $\alpha < \kappa$---i.e.\ if $V_\alpha \cap \mathbb{P}_\alpha$ is a regular suborder of $\mathbb{P}_\alpha$---the $V_\alpha \cap \mathbb{P}_\alpha$-name that we force with is the following poset, assuming that $\alpha = \mu^+$ in $V^{V_\alpha \cap \mathbb{P}_\alpha}$:
\[
\dot{Q}_\alpha:= \Big( \mathbb{C}_{\text{Silv}}(\alpha, < \kappa) \ * \  \dot{\mathbb{S}}^{H_{\alpha^{++}}}_{r_\alpha,\phi}  \Big)^{V[\dot{g}_\alpha]}
\]
where $\dot{g}_\alpha$ is the $V_\alpha \cap \mathbb{P}_\alpha$-name for its generic object, $\mathbb{C}_{\text{Silv}}(\alpha, < \kappa)$ is the Silver collapse (see Section \ref{sec_Examples}) that turns $\kappa$ into $\alpha^+$, and $r_\alpha$ is some $(V_\alpha \cap \mathbb{P}_\alpha)*\dot{\mathbb{C}}_{\text{Silv}}(\alpha, < \kappa)$-name for a parameter which witnesses the analogue of \eqref{eq_DefinablePlus2} from the statement of Theorem \ref{thm_IdealAbsorb}.  Note that, from the point of view of the model $V[g_\alpha]^{\mathbb{C}_{\text{Silv}}(\alpha, < \kappa)}$, the poset $\dot{\mathbb{S}}^{H_{\alpha^{++}}}_{r_\alpha,\phi}$ satisfies the assumptions of Theorem \ref{thm_IdealAbsorb}.

It is routine to verify that  
\begin{equation*}
\Vdash_{V_\alpha \cap \mathbb{P}_\alpha} \ \dot{Q}_\alpha \text{ is } \kappa \text{-cc} \text{, } <\mu \text{-closed, } \text{ and has size } \kappa.
\end{equation*}
Then by Theorem \ref{thm_WeakCompact_Kunen_MainThm_Precise} it follows that
\begin{equation}\label{eq_P_is_KappaCC}
\mathbb{P}:= \mathbb{P}_\kappa \text{ is } \kappa \text{-cc}
\end{equation} 
and Lemma \ref{lem_LessMuClosedIteration} ensures that $\mathbb{P}$ is $<\mu$ closed.  Because $\mathbb{P}_0$ collapses all cardinals between $\mu$ and $\kappa$, these facts together imply
\begin{equation}\label{eq_KappaBecomesMuPlus}
\Vdash_{\mathbb{P}} \ \kappa = \mu^+.
\end{equation}

For reasons that will be discussed in Section \ref{sec_HugeVersusAH}, we will primarily work with the \emph{almost huge embedding $j_{\vec{U}}: V \to N_{\vec{U}}$ derived from $j$}, as in Fact \ref{fact_DerivedAlmostHuge}, rather than the huge embedding $j$ (though $j$ will play a key role).  Let $k: N_{\vec{U}} \to N$ be the map from Fact \ref{fact_DerivedAlmostHuge}, and recall from Fact \ref{fact_DerivedAlmostHuge} that $\text{crit}(k) = (\delta^+)^{N_{\vec{U}}}$ (recall $\delta = j(\kappa)$).

By \eqref{eq_P_is_KappaCC} and Lemma \ref{lem_SubsetVkappa} it follows that $\kappa$ is an active stage of $j_{\vec{U}}(\vec{\mathbb{P}})$.  Let $G$ be $(V,\mathbb{P})$-generic.  Then by the definition of the iteration, combined with \eqref{eq_KappaBecomesMuPlus}, stage $\kappa$ of the $j_{\vec{U}}(\vec{\mathbb{P}})$ iteration forces with the following poset:
\begin{equation}\label{eq_2_step_in_N}
\Big(\mathbb{C}_{\text{Silv}}(\kappa, < \delta)  \ * \  \dot{\mathbb{S}}^{H_{\kappa^{++}}}_{j_{\vec{U}}(\vec{r})(\kappa), \phi} \Big)^{N_{\vec{U}}[G]}.
\end{equation}

Because $N_{\vec{U}}$ is $<\delta$ closed in $V$, and $G$ is generic for a $\delta$-cc poset, then $N_{\vec{U}}[G]$ is $<\delta$ closed in $V[G]$.  In particular, $V[G]$ and $N_{\vec{U}}[G]$ have the same $\kappa$ sequences, and hence $\mathbb{C}_{\text{Silv}}(\kappa, < \delta)$ is computed the same in $V[G]$ and $N_{\vec{U}}[G]$.  Let $H$ be $(V[G], \mathbb{C}_{\text{Silv}}(\kappa, < \delta))$-generic, and consider the evaluation of the 2nd step of \eqref{eq_2_step_in_N}; i.e.\ the poset
\begin{equation*}
\Big( \mathbb{S}^{H_{\kappa^{++}}}_{r, \phi} \Big)^{N_{\vec{U}}[G*H]}
\end{equation*}
where $r$ is the evaluation of $j_{\vec{U}}(\vec{r})(\kappa)$ by $G*H$; so $r$ is a parameter in $(H_{\kappa^{++}})^{N_{\vec{U}}[G*H]} = (H_{\delta^+})^{N_{\vec{U}}[G*H]}$.  Now we recall the assumptions about the parameter $r$, namely
\begin{align}\label{eq_RecallWhatNBelieves}
\begin{split}
N_{\vec{U}}[G*H] \models  \mathbb{S}^{H_{\kappa^{++}}}_{r,\phi}= \big\{ z \in H_{\kappa^{++}} \ |  \ (H_{\kappa^{++}}, \in) \models \phi(z,r) \big\}  \\
\text{ is a } <\mu \text{-closed, } \kappa^+ \text{-cc poset} \text{ of size at most } \kappa^+.
\end{split}
\end{align}

This implies that $(\mathbb{S}^{H_{\kappa^{++}}}_{r,\phi})^{N_{\vec{U}}[G*H]}$ is actually an element of $(H_{\kappa^{++}})^{N_{\vec{U}}[G*H]}$ and that
\begin{align}\label{eq_LocallyBelieved}
\begin{split}
(H_{\kappa^{++}})^{N_{\vec{U}}[G*H]} \models \ \ (\mathbb{S}^{H_{\kappa^{++}}}_{r,\phi})^{N_{\vec{U}}[G*H]} \text{ is a } <\mu \text{-closed, } \\
\kappa^+ \text{-cc poset of size at most } \kappa^+.
\end{split}
\end{align}

Let $\mathbb{S}:= (\mathbb{S}^{H_{\kappa^{++}}}_{r,\phi})^{N_{\vec{U}}[G*H]}$.  By the definition of $j_{\vec{U}}(\vec{\mathbb{P}})$, there is an $\iota \in N_{\vec{U}}$ such that
\[
\iota: \mathbb{P}*\dot{\mathbb{C}}_{\text{Silv}}(\kappa, < \delta)*\dot{\mathbb{S}} \to j_{\vec{U}}(\mathbb{P}) \text{ is a regular embedding and is the identity on } \mathbb{P}.
\]

By Theorem \ref{thm_Duality_AH}, in $V[G*H]$ there is a normal ideal $\mathcal{I}(j_{\vec{U}})$ on $\kappa$ such that
\begin{equation}\label{eq_ForcingEquiv}
\wp(\kappa)/\mathcal{I}(j_{\vec{U}}) \text{ is forcing equivalent to } j_{\vec{U}}(\mathbb{P})/ \iota[G*H].
\end{equation}
By standard arguments involving quotient forcings, the map $\iota$ can be used to show that the poset $\mathbb{S} = \dot{\mathbb{S}}^{V[G*H]}$ is absorbed as a regular suborder of $j_{\vec{U}}(\mathbb{P})/ \iota[G*H]$.  Combining this with \eqref{eq_ForcingEquiv} yields:
\begin{equation}\label{eq_RegEmbeddingThere}
V[G*H] \models \ \ \text{There is a regular embedding } e: \mathbb{S} \to \text{ro} \big( \wp(\kappa)/\mathcal{I}(j_{\vec{U}}) \big)
\end{equation}
where $\text{ro}(-)$ denotes the boolean completion.

\begin{remark}
Regarding the ``ro" appearing in \eqref{eq_RegEmbeddingThere}:  below we will show that in fact $\mathcal{I}(j_{\vec{U}})$ is saturated, which implies by Fact \ref{fact_SaturatedImpliesCompleteBA} that $\wp(\kappa)/\mathcal{I}(j_{\vec{U}})$ is a complete boolean algebra, i.e.\ that $\wp(\kappa)/\mathcal{I}(j_{\vec{U}})$ is isomorphic to (not just densely embeddable into) the complete boolean algebra $\text{ro} \big( \wp(\kappa)/\mathcal{I}(j_{\vec{U}}) \big)$.  So the ``ro" in \eqref{eq_RegEmbeddingThere} will ultimately be redundant.
\end{remark}

\begin{remark}
Using the Duality Theorey from Foreman~\cite{MattHandbook}, it can be arranged that the regular embedding $e$ from \eqref{eq_RegEmbeddingThere} has the following additional property:  that whenever $K$ is $\big( V[G*H],\wp(\kappa)/\mathcal{I}(j_{\vec{U}}) \big)$-generic and $j_K: V[G*H] \to \text{ult}(V[G*H],K)$ is the generic ultrapower map, then $e^{-1}[K]$ is an element of $\text{ult}(V[G*H],K)$ (although $K$ will never be an element of $\text{ult}(V[G*H],K)$).   However we choose not to include this additional feature in the theorem, in order to avoid lengthy detours into Duality Theory.
\end{remark}

Now it remains to prove the following two claims:
\begin{globalClaim}\label{clm_Saturated}
The ideal $\mathcal{I}(j_{\vec{U}})$ is saturated, i.e.
\[
V[G*H] \models \wp(\kappa)/\mathcal{I}(j_{\vec{U}}) \text{ is } \delta = \kappa^+ \text{-cc}.
\]  
\end{globalClaim}

\begin{globalClaim}\label{clm_EqualityToProve} 
The poset $\mathbb{S}$---which, recall, is the poset  $(\mathbb{S}^{H_{\kappa^{++}}}_{r,\phi})^{N_{\vec{U}}[G*H]}$ which was defined in $N_{\vec{U}}[G*H]$---has the required definability properties from the point of view of $V[G*H]$.  In other words, it remains to prove that the following holds from the point of view of $V[G*H]$:
\begin{align*}
  \mathbb{S} = \{ z \in (H_{\kappa^{++}})^{V[G*H]} \ | \   \big(  (H_{\kappa^{++}})^{V[G*H]}, \in \big) \models \phi(r,z)   \} \\
\textbf{and } \mathbb{S} \text{ is a } <\mu \text{-closed, } \kappa^+ \text{-cc poset of size at most } \kappa^+.
\end{align*}

\end{globalClaim}

\begin{remark}\label{rem_RegardingHuge}
The equality 
\[
\mathbb{S} = \{ z \in (H_{\kappa^{++}})^{V[G*H]} \ | \   \big(  (H_{\kappa^{++}})^{V[G*H]}, \in \big) \models \phi(r,z)   \}
\]
from Claim \ref{clm_EqualityToProve} does \textbf{not} follow automatically from \eqref{eq_RecallWhatNBelieves}, because 
\[
(\kappa^{++})^{V[G*H]} = \delta^{+V} > \delta^{+N_{\vec{U}}} =  (\kappa^{++})^{N_{\vec{U}}[G*H]}
\]
where the inequality is by Fact \ref{fact_DerivedAlmostHuge}.  We seem to need some additional assumption, such as the huge embedding used in the argument below.
\end{remark}

We will use the huge embedding $j$ to help prove Claims \ref{clm_Saturated} and \ref{clm_EqualityToProve}.  Let 
$k: N_{\vec{U}} \to N$ be the map from Fact \ref{fact_DerivedAlmostHuge}.  Recall that
\begin{equation}\label{eq_CritOfk}
\tau:= \text{crit}(k) = (\delta^+)^{N_{\vec{U}}} = (\kappa^{++})^{N_{\vec{U}}[G*H]}.
\end{equation}

\subsubsection{Proof of Claim \ref{clm_Saturated}}

By \eqref{eq_ForcingEquiv}, proving Claim \ref{clm_Saturated} is equivalent to showing that $V[G*H] \models$ ``the poset $j_{\vec{U}}(\mathbb{P})/\iota[G*H]$ is $\delta$-cc".  By Fact \ref{fact_CannotAddAntichains}, in turn it suffices to prove that $V \models$ ``the poset $j_{\vec{U}}(\mathbb{P})$ is $\delta$-cc".  But because $j_{\vec{U}}(\mathbb{P}) \subset j_{\vec{U}}(V_\kappa) = (V_\delta)^{N_{\vec{U}}}$, then by \eqref{eq_CritOfk} the map $k: N_{\vec{U}} \to N$ fixes $j_{\vec{U}}(\mathbb{P})$.  Also because $\mathbb{P}$ is $\kappa$-cc in $V$ then $j_{\vec{U}}(\mathbb{P})$ is $\delta$-cc in $N_{\vec{U}}$.  Then by elementarity of $k$, $k(j_{\vec{U}}(\mathbb{P})) = j_{\vec{U}}(\mathbb{P})$ is $k(\delta) = \delta$-cc from the point of view of $N$.  Because $N$ is closed under $\delta$ sequences in $V$, it follows that $j_{\vec{U}}(\mathbb{P})$ is really $\delta$-cc in $V$ as well, which completes the proof.

\subsubsection{Proof of Claim \ref{clm_EqualityToProve}}

First note that $k$ lifts to an elementary $k: N_{\vec{U}}[G*H] \to N[G*H]$ because $G*H$ is generic for a poset which is a subset of $V_\delta^{N_{\vec{U}}}$ and $\delta < \text{crit}(k)$.  Because $\mathbb{S} \in (H_\tau)^{N_{\vec{U}}[G*H]}$ and $\text{crit}(k) = \tau$ then $k(\mathbb{S}) = \mathbb{S}$; similarly because $r \in (H_\tau)^{N_{\vec{U}}[G*H]}$ then $k(r) = r$.  Then by \eqref{eq_RecallWhatNBelieves} and elementarity of $k$ it follows that
\begin{align}\label{eq_HugeBelieves}
\begin{split}
N[G*H] \models k(\mathbb{S}) = \mathbb{S} = \{ z \in \big( H_{\kappa^{++}}\big)^{N[G*H]} \ | \  \Big( \big( H_{\kappa^{++}}\big)^{N[G*H]}, \in \Big) \models \phi(z,r)   \} \\
\text{ and this poset is } <\mu \text{-closed, } \delta\text{-cc, and of size } \le \delta.
\end{split}
\end{align}

Because the $\models$ relation is absolute, \eqref{eq_HugeBelieves} implies
\begin{equation}\label{eq_V_Believes}
V[G*H] \models \mathbb{S} = \{ z \in \big( H_{\kappa^{++}}\big)^{N[G*H]} \ | \  \Big( \big( H_{\kappa^{++}}\big)^{N[G*H]}, \in \Big) \models \phi(z,r)   \}.
\end{equation}

Because $N$ is closed under $\delta$ sequences and $G*H$ is generic for a $\delta$-cc poset, $N[G*H]$ is closed under $\delta$ sequences in $V[G*H]$; in particular $\delta^{+V[G*H]} = \delta^{+N[G*H]}$ and
\begin{equation}
(H_{\kappa^{++}})^{V[G*H]} = (H_{\delta^+})^{V[G*H]} =(H_{\delta^+})^{N[G*H]}  =  (H_{\kappa^{++}})^{N[G*H]}.
\end{equation}
Thus by substituting into \eqref{eq_V_Believes} we obtain
\begin{equation}
V[G*H] \models \mathbb{S} = \{ z \in \big( H_{\kappa^{++}}\big)^{V[G*H]} \ | \  \Big( \big( H_{\kappa^{++}}\big)^{V[G*H]}, \in \Big) \models \phi(z,r)   \}.
\end{equation}
This shows that $\mathbb{S}$ indeed has the correct definition from the point of view of $V[G*H]$.  Finally, because $N[G*H]$ is $\delta$-closed in $V[G*H]$, the fact that $N[G*H]$ believes $\mathbb{S}$ is $<\mu$-closed, $\delta$-cc, and of size $\le \delta$ is upward absolute to $V[G*H]$.

\subsubsection{A remark about the use of the huge embedding}\label{sec_HugeVersusAH}

A remark is in order regarding why we insisted on using the almost huge embedding $j_{\vec{U}}$ rather than the huge embedding $j$ to construct the saturated ideal $\mathcal{I}(j_{\vec{U}})$ in the proof of Theorem \ref{thm_IdealAbsorb}.  It certainly is possible to use liftings of the huge embedding $j$, together with a ``pseudo-generic tower" construction, to define a saturated ideal $\mathcal{I}(j)$ on $\kappa$; this was Kunen's original argument.  However, that argument relied on the coarse fact that there exists an incompatibility-preserving embedding $i$ from $\wp(\kappa)/\mathcal{I}(j)$ into the $\delta=\kappa^+$-cc poset 
\[
\mathbb{P}^*:= j(\mathbb{P})/\iota[G*H] = j_{\vec{U}}(\mathbb{P})/\iota[G*H]
\]
which guarantees that $\wp(\kappa)/\mathcal{I}(j)$ is $\kappa^+$-cc, i.e.\ that $\mathcal{I}(j)$ is saturated.  However, it is not clear if the map $i$ is a dense embedding, or even a regular embedding.  So in particular it is not clear if $\mathbb{S}$---which is a regular suborder of $\mathbb{P}^*$---can be regularly embedded into $\wp(\kappa)/\mathcal{I}(j)$.  On the other hand, $\wp(\kappa)/\mathcal{I}(j_{\vec{U}})$ \textbf{is} forcing equivalent to $\mathbb{P}^*$ and thus absorbs $\mathbb{S}$ as a regular suborder, which is why we worked with $\mathcal{I}(j_{\vec{U}})$ instead of $\mathcal{I}(j)$.  Nonetheless, the map $j$ played an important background role in the proofs of Claims \ref{clm_Saturated} and \ref{clm_EqualityToProve}.

\subsection{Proof of Theorem \ref{thm_ChangTree}}

Now we proceed to the proof of Theorem \ref{thm_ChangTree}.  For regular cardinals $\kappa < \delta$, let $\text{Sacks}(\kappa, \delta)$ denote the $\le \kappa$-support iteration of length $\delta$ which uses $\text{Sacks}(\kappa)$ at each step (see Definition \ref{def_Sacks} for the definition of $\text{Sacks}(\kappa)$).  We use the following theorem of Kanamori, which generalizes a theorem from Baumgartner-Laver~\cite{MR556894}: 
\begin{theorem}[Kanamori~\cite{MR593029}]\label{thm_Kanamori}
Suppose $\kappa < \delta$ are regular uncountable, $\delta$ is weakly compact, and $\Diamond_\kappa$ holds.  Then the $\le \kappa$ support iteration of length $\delta$ of $\text{Sacks}(\kappa)$ is $\delta$-cc, and forces the following:
\begin{itemize}
 \item $\delta = \kappa^{++}$;
 \item The Tree property holds at $\delta = \kappa^{++}$.
\end{itemize}
\end{theorem}

For regular $\kappa < \delta$ let $\text{Sacks}_{\Diamond}(\kappa, \delta)$ denote the $\le \kappa$-support iteration 
\[
\langle \mathbb{R}_\gamma, \dot{\mathbb{Q}}_\beta \ | \ \gamma \le \delta, \beta < \delta \rangle
\]
where $\mathbb{Q}_0 = \text{Add}(\kappa)$ and $\dot{\mathbb{Q}}_\gamma = \text{Sacks}^{\mathbb{R}_\gamma}(\kappa)$ for all $\gamma \in (0,\delta)$.  It is a standard fact that $\text{Add}(\kappa)$ forces $\Diamond(\kappa)$; so by Theorem \ref{thm_Kanamori} it follows that:
\begin{corollary}\label{cor_ForcesTP}
If $\kappa < \delta$ are regular uncountable and $\delta$ is weakly compact, then $\text{Sacks}_{\Diamond}(\kappa, \delta)$ is $\delta$-cc and forces that $\delta = \kappa^{++}$ and the Tree Property holds at $\kappa^{++}$.
\end{corollary}

Assume $j: V \to N$ is a huge embedding with $\text{crit}(j) = \kappa$ and $j(\kappa) = \delta$.  Fix a regular $\mu < \kappa$.  Define a $<\mu$-support universal Kunen iteration  
\[
\langle \mathbb{P}_\alpha, \dot{Q}_\beta \ | \ \alpha \le \kappa, \beta < \kappa \rangle
\]
(as in Definition \ref{def_KunenIteration}) as follows:  $\mathbb{P}_0$ is $\text{Col}(\mu, < \kappa)$.  For active $\alpha < \kappa$---i.e.\ if $V_\alpha \cap \mathbb{P}_\alpha$ is a regular suborder of $\mathbb{P}_\alpha$---the $V_\alpha \cap \mathbb{P}_\alpha$-name that we force with is the following poset:
\[
\dot{Q}_\alpha:= \Big(  \text{Sacks}_\Diamond(\alpha, \kappa \Big)^{V[\dot{g}_\alpha]}
\]
where $\dot{g}_\alpha$ is the $V_\alpha \cap \mathbb{P}_\alpha$-name for its generic object.  Note that if $\alpha < \kappa$ and $g_\alpha$ is generic for $V_\alpha \cap \mathbb{P}_\alpha$, then $\kappa$ is still weakly compact in $V[g_\alpha]$ (because $g_\alpha$ is generic for a $<\kappa$-sized poset), and so by Corollary \ref{cor_ForcesTP} it follows that:
\begin{equation*}
\Vdash_{V_\alpha \cap \mathbb{P}_\alpha} \ \dot{Q}_\alpha \text{ is } \kappa \text{-cc and forces } \kappa = \alpha^{++} \ \wedge \ \text{TP}(\kappa) .
\end{equation*}
By Theorem \ref{thm_WeakCompact_Kunen_MainThm_Precise}, $\mathbb{P}:= \mathbb{P}_\kappa$ is $\kappa$-cc; then by Lemma \ref{lem_SubsetVkappa} it follows that $\kappa$ is an active stage of $j(\vec{\mathbb{P}})$.  Moreover because $\mathbb{P}_0$ collapses all cardinals in the interval $(\mu,\kappa)$ and each $\dot{Q}_\alpha$ is forced to be $<\mu$-closed, then $\kappa = \mu^+$ in $V^{\mathbb{P}}$.  

Let $G$ be $(V,\mathbb{P})$-generic.  Then stage $\kappa$ of the $j(\vec{\mathbb{P}})$ iteration forces with the following poset:
\begin{equation}
\Big(\text{Sacks}_\Diamond(\kappa, \delta)\Big)^{N[G]}.
\end{equation}
So there is an $e \in N[G]$ such that
\[
e: \Big( \text{Sacks}_\Diamond(\kappa, \delta)\Big)^{N[G]} \to j(\mathbb{P})/G \text{ is a regular embedding}.
\]

Because $N$ is closed under $\delta$ sequences in $V$ (really all we need for this part is that it is closed under $<\delta$ sequences) it follows that 
\begin{equation}\label{eq_SamePoset}
\mathbb{Q}:=\Big(\text{Sacks}_\Diamond(\kappa, \delta)\Big)^{V[G]} = \Big(\text{Sacks}_\Diamond(\kappa, \delta)\Big)^{N[G]}.
\end{equation}

Note that, from the point of view of $V[G]$, the poset $\mathbb{Q}$ is a $\delta$-length, $\le \kappa$ support iteration of $<\kappa$-directed closed posets, where the poset at each step has size $<\delta$.  That each step is $<\kappa$-directed closed follows from Lemma \ref{lem_SacksDirectedClosed}.

Let $K$ be generic over $V[G]$ for the poset $\mathbb{Q}$.  Then by Corollary \ref{cor_ForcesTP}, in $V[G*K]$ the Tree Property holds at $\kappa^{++}$.  Also $\kappa = (\mu^+)^{V[G*K]}$, because $\mathbb{P}$ forces $\kappa = \mu^+$ and $\mathbb{Q}$ is $<\kappa$ (directed) closed.  Also $\mu$ is still a cardinal by Lemma \ref{lem_LessMuClosedIteration}.  It remains to prove that 
\[
V[G*K] \models (\kappa^{++}, \kappa) \twoheadrightarrow (\kappa, \mu).
\]
Let $G'$ be generic over $V[G*K]$ for 
\[
(j(\mathbb{P})/G)/e[K].
\]
Then in $V[G']$, the map $j$ lifts to
\[
j: V[G] \to N[G'].
\]
We want to further lift $j$ to domain $V[G*K]$ in some forcing extension.  Note in particular that $K \in N[G']$, because $e \in N[G]$ and $K = e^{-1}[G']$.  Since $N$ is closed under $\delta$ sequences and $j(\mathbb{P})$ is $\delta$-cc, it follows that $N[G']$ is closed under $\delta$ sequences from $V[G']$; in particular $j \restriction A \in N[G']$ for every $A \in V[G]$ of size $\delta$.  So the map $j: V[G] \to N[G']$ and the iteration $\mathbb{Q}$ satisfy the requirements of the following lemma, with $V[G]$ playing the role of $W$ and $N[G']$ playing the role of $W'$.

\begin{lemma}[implicit in Kunen~\cite{MR495118}]
Suppose $j:W \to W'$ is a (possibly external) elementary embedding with critical point $\kappa$; let $\delta:= j(\kappa)$ and assume that $\delta$ is inaccessible in $W$.  Suppose $\vec{\mathbb{Q}} = \langle \mathbb{Q}_\gamma, \dot{R}_\beta \ | \ \gamma \le \delta, \beta < \delta \rangle$ is a $\le \kappa$-support iteration in $W$ such that
\begin{equation}
\forall \gamma < \delta \ \Vdash^W_{\mathbb{Q}_\gamma} \dot{R}_\gamma \text{ is } < \kappa \text{ directed closed and of size } < \delta.
\end{equation}

Assume that the iteration $\vec{\mathbb{Q}}$ is also an element of $W'$, that  
\begin{equation}
\forall A \in \wp^W(V_\delta) \ \ \  j \restriction A \in W'
\end{equation}
and that there exists some $K \in W'$ which is $(W,\mathbb{Q}_\delta)$-generic.

Then:
\begin{itemize}
 \item there is a condition $m_K \in j(\mathbb{Q}_\delta)$ such that if $K'$ is $(W',j(\mathbb{Q}_\delta))$-generic and $m_K \in K'$, then $j[K] \subset K'$ and thus $j$ can be lifted to $\tilde{j}: W[K] \to W'[K']$.
 \item $W[K] \models (\delta, \kappa) \twoheadrightarrow (\kappa, < \kappa)$
\end{itemize} 
\end{lemma}
\begin{proof}
The master condition $m_K$ will essentially be the condition with support $j[\delta]$ such that for each $\gamma < \delta$, the $j(\gamma)$-th component of $m_K$ is a lower bound for the pointwise image of the $\gamma$-th component of $K$; this is possible because each component of the iteration is of size $<\delta$ and the image of the components of the iteration are $<\delta$-directed closed.

More precisely, we will recursively define a sequence
\[
\vec{m} = \langle m_\gamma \ | \ \gamma \le \delta \rangle
\]
 and inductively verify that:
\begin{enumerate}
 \item $m_\gamma \in j(\mathbb{Q}_\gamma)$ and has support exactly $j[\gamma]$;
 \item Whenever $\gamma < \gamma'$ then $m_\gamma = m_{\gamma'} \restriction j(\gamma)$;
 \item If $K'_{j(\gamma)}$ is generic for $j(\mathbb{Q}_\gamma)$ over $W'$ and $m_\gamma \in K'_{j(\gamma)}$, then $j[K|\gamma] \subset K'_{j(\gamma)}$ and thus $j$ lifts to an embedding $j: W[K|\gamma] \to W'[K'_{j(\gamma)}]$.
\end{enumerate}

The recursion goes as follows.  If $\gamma$ is a limit ordinal then $m_\gamma$ is simply the union of the $m_\beta$ for $\beta < \gamma$.  If $\gamma$ is a successor ordinal, say $\gamma = \bar{\gamma}+1$, then we define $m_{\bar{\gamma}+1}$ as $m_{\bar{\gamma}} \cup \{ \big(j(\bar{\gamma}), \dot{k}_{\bar{\gamma}} \big) \}$ where $\dot{k}_{\bar{\gamma}}$ is the $j(\mathbb{Q}_{\bar{\gamma}})$-name defined conditionally as follows (here $\dot{K}'_{j(\bar{\gamma})}$ denotes the $j(\mathbb{Q}_{\bar{\gamma})}$-name for its generic object):
\begin{itemize}
 \item If $m_{\bar{\gamma}} \notin \dot{K}'_{j(\bar{\gamma})}$ then $\dot{k}_{\bar{\gamma}}$ is the trivial condition of $j(\dot{R}_{\bar{\gamma}})$;
 \item If $m_{\bar{\gamma}} \in \dot{K}'_{j(\bar{\gamma})}$ then let 
\[ 
 j_{\bar{\gamma}}: W[K|\bar{\gamma}] \to  W'[\dot{K}'_{j(\bar{\gamma})}]
\]
denote the lifting which is guaranteed by the induction hypothesis.  Let $H_{\bar{\gamma}}$ be the $\bar{\gamma}$-th component of $K$, i.e.\ 
 \[
 H_{\bar{\gamma}}:= \{ \big(r(\bar{\gamma})\big)_{K|\bar{\gamma}} \ | \  r \in K|(\bar{\gamma}+1)   \} \subset \big( \dot{R}_{\bar{\gamma}} \big)_{K|\bar{\gamma}}.
 \]
Notice that $H_{\bar{\gamma}} \in W'$ because $K \in W'$.  Because $\big( \dot{R}_{\bar{\gamma}} \big)_{K|\bar{\gamma}}$ has size $<\delta$ and $j_{\bar{\gamma}}$ is elementary, then $j_{\bar{\gamma}}[H_{\bar{\gamma}}]$ is a directed set of conditions of size $<\delta$ in the model $W'[\dot{K}'_{j(\bar{\gamma})}]$.  And because $\big( \dot{R}_{\bar{\gamma}} \big)_{K|\bar{\gamma}}$ is $<\delta$-directed closed in $W'[\dot{K}'_{j(\bar{\gamma})}]$ then $j_{\bar{\gamma}}[H_{\bar{\gamma}}]$ has a lower bound; let $\dot{k}_{\bar{\gamma}}$ be any such lower bound for $j_{\bar{\gamma}}[H_{\bar{\gamma}}]$.
\end{itemize}

It is now routine to verify the inductive assumptions.  Note that $m_\delta$ has support $j[\delta]$, which is acceptable because $j(\vec{\mathbb{Q}})$ is a $\le \delta$ support iteration and $|j[\delta]|=\delta$.

Finally we prove that $W[K] \models (\delta, \kappa) \twoheadrightarrow (\kappa, < \kappa)$.  In $W[K]$ fix a structure $\mathfrak{A} = (\delta, \dots)$ in a countable language.  Let $K'$ be $\big(W', j(\mathbb{Q}_\delta) \big)$-generic with $m_K \in K'$ and let $j_\delta: W[K] \to W'[K']$ be the lifting of $j$.  Elementarity of $j_\delta$ ensures that $j_\delta[\delta] = j[\delta] \prec j_\delta(\mathfrak{A})$.  Moreover $j[\delta] \in W'$ and $W'[K']$ believes that $j[\delta] \cap j(\kappa) = \kappa \in j(\kappa)$ and that $|j[\delta]| = \delta = j(\kappa)$; so
\[
W'[K'] \models \ \exists X \prec j_\delta(\mathfrak{A}) \ \ \ X \cap j(\kappa) \in j(\kappa) \text{ and } |X|  = j(\kappa).
\]
So by elementarity of $j_\delta$:
\[
W[K] \models \ \exists X \prec \mathfrak{A} \ \ \ X \cap \kappa \in \kappa \text{ and } |X|  = \kappa
\]
which completes the proof.
\end{proof}

\begin{remark}
The model just constructed actually has a $\kappa^{++}$-saturated ideal on 
\[
\{ z \in [\kappa^{++}]^\kappa \ | \ z \cap \kappa \in \kappa   \}.
\]
This can be proved using the ``pseudogeneric tower" argument from Kunen~\cite{MR495118} which is described in detail in Foreman~\cite{MattHandbook}.
\end{remark}

\section{Questions}

\begin{question}
Suppose we alter the assumptions of Theorem \ref{thm_WeakCompact_Kunen_MainThm_Precise} as follows:
\begin{itemize}
 \item Weaken the large cardinal assumption (e.g.\ only assume that $\kappa$ is Mahlo)
 \item Strengthen the assumptions on $\dot{Q}_\alpha$ by requiring them to be $\kappa$-Knaster instead of merely $\kappa$-cc.  More precisely, assume that for all active $\alpha < \kappa$ that $\Vdash_{V_\alpha \cap \mathbb{P}_\alpha}$ $\dot{Q}_\alpha$ is $\kappa$-Knaster.  
\end{itemize}
Can we conclude that $\mathbb{P}_\kappa$ has the $\kappa$-cc?  Note that if $\vec{\mathbb{P}}$ is a finite support iteration, then the answer is yes.  
\end{question}

Remark \ref{rem_RegardingHuge} in the proof of Theorem \ref{thm_IdealAbsorb} raises the following question:
\begin{question}
Is Theorem \ref{thm_IdealAbsorb} still true if $\kappa$ is merely assumed to be \textbf{almost} huge in $V$?
\end{question}

\end{document}